\documentclass[letter, reqno,11pt]{amsart}

\usepackage[usenames,dvipsnames]{color}
\usepackage[utf8]{inputenc}
\usepackage{amsthm,amsfonts,amssymb,amsmath,amsxtra}
\usepackage[all]{xy}
\SelectTips{cm}{}

\usepackage{nccrules}

\usepackage{verbatim}
\usepackage{xr-hyper}
\usepackage{varioref}
\usepackage{hyperref}
\usepackage[nameinlink]{cleveref}
\usepackage{tikz-cd}
\usepackage[margin=1.25in]{geometry}
\usepackage{mathrsfs}
\usepackage{xcolor}

\usepackage[shortlabels]{enumitem}
\usepackage{lipsum}
\usetikzlibrary{positioning,chains,fit,shapes,calc,arrows,patterns,external,shapes.callouts,graphs}

\usepackage[style=alphabetic]{biblatex}
\RequirePackage{xspace}
\RequirePackage{etoolbox}
\RequirePackage{varwidth}
\RequirePackage{enumitem}
\RequirePackage{tensor}
\RequirePackage{mathtools}
\RequirePackage{longtable}
\RequirePackage{multirow}

\def\ge{\geqslant}
\def\le{\leqslant}
\def\a{\alpha}

\def\s{\sigma}

\def\i{^{-1}}

\def\<{\langle}
\def\>{\rangle}

\newcommand{\bG}{\mathbf G}

\newcommand{\BA}{\ensuremath{\mathbb {A}}\xspace}

\newcommand{\BC}{\ensuremath{\mathbb {C}}\xspace}

\newcommand{{\BG}}{\ensuremath{\mathbb {G}}\xspace}

\newcommand{{\BK}}{\ensuremath{\mathbb {K}}\xspace}

\newcommand{\BS}{\ensuremath{\mathbb {S}}\xspace}

\newcommand{\BZ}{\ensuremath{\mathbb {Z}}\xspace}

\newcommand{\CC}{\ensuremath{\mathcal {C}}\xspace}

\newcommand{\CI}{\ensuremath{\mathcal {I}}\xspace}

\newcommand{\CO}{\ensuremath{\mathcal {O}}\xspace}
\newcommand{\CP}{\ensuremath{\mathcal {P}}\xspace}

\newcommand{\CS}{\ensuremath{\mathcal {S}}\xspace}

\newcommand{\CU}{\ensuremath{\mathcal {U}}\xspace}

\newcommand{\tW}{\widetilde{W}}

\newcommand{\dimm}{\dim}

\DeclareMathOperator{\supp}{supp}

\newcommand{\wt}{\text{wt}}

\def\tW{\tilde W}

\def\kk{\mathbf k}
\newcommand{\ii}{\mathbf i}

\newcommand{\Z}{\mathbb Z}
\newcommand{\R}{\mathbb R}

\newcommand{\ee}{\widetilde e}
\newcommand{\ff}{\widetilde f}

\newcommand{\inv}{{\mathrm{inv}}}

\def\kk{\mathbf k}

\newtheorem{theorem}{Theorem}
\newtheorem{proposition}[theorem]{Proposition}
\newtheorem{lemma}[theorem]{Lemma}

\newtheorem{corollary}[theorem]{Corollary}

\theoremstyle{definition}

\newtheorem{definition}[theorem]{Definition}
\newtheorem{example}[theorem]{Example}

\newtheorem{remark}[theorem]{Remark}

\numberwithin{equation}{section}
\numberwithin{theorem}{section}

\setitemize[0]{leftmargin=*,itemsep=\the\smallskipamount}
\setenumerate[0]{leftmargin=*,itemsep=\the\smallskipamount}

\renewcommand{\to}{%
   \ifbool{@display}{\longrightarrow}{\rightarrow}%
   }
\let\shortmapsto\mapsto
\renewcommand{\mapsto}{%
   \ifbool{@display}{\longmapsto}{\shortmapsto}%
   }
\newlength{\olen}
\newlength{\ulen}
\newlength{\xlen}
\newcommand{\xra}[2][]{%
   \ifbool{@display}%
      {\settowidth{\olen}{$\overset{#2}{\longrightarrow}$}%
       \settowidth{\ulen}{$\underset{#1}{\longrightarrow}$}%
       \settowidth{\xlen}{$\xrightarrow[#1]{#2}$}%
       \ifdimgreater{\olen}{\xlen}%
          {\underset{#1}{\overset{#2}{\longrightarrow}}}%
          {\ifdimgreater{\ulen}{\xlen}%
             {\underset{#1}{\overset{#2}{\longrightarrow}}}
             {\xrightarrow[#1]{#2}}}}%
      {\xrightarrow[#1]{#2}}
   }
\makeatother
\newcommand{\xyra}[2][]{%
   \settowidth{\xlen}{$\xrightarrow[#1]{#2}$}%
   \ifbool{@display}%
      {\settowidth{\olen}{$\overset{#2}{\longrightarrow}$}%
       \settowidth{\ulen}{$\underset{#1}{\longrightarrow}$}%
       \ifdimgreater{\olen}{\xlen}%
          {\mathrel{\xymatrix@M=.12ex@C=3.2ex{\ar[r]^-{#2}_-{#1} &}}}%
          {\ifdimgreater{\ulen}{\xlen}%
             {\mathrel{\xymatrix@M=.12ex@C=3.2ex{\ar[r]^-{#2}_-{#1} &}}}
             {\mathrel{\xymatrix@M=.12ex@C=\the\xlen{\ar[r]^-{#2}_-{#1} &}}}}}%
      {\mathrel{\xymatrix@M=.12ex@C=\the\xlen{\ar[r]^-{#2}_-{#1} &}}}%
   }
\makeatletter
\newcommand{\xla}[2][]{%
   \ifbool{@display}%
      {\settowidth{\olen}{$\overset{#2}{\longleftarrow}$}%
       \settowidth{\ulen}{$\underset{#1}{\longleftarrow}$}%
       \settowidth{\xlen}{$\xleftarrow[#1]{#2}$}%
       \ifdimgreater{\olen}{\xlen}%
          {\underset{#1}{\overset{#2}{\longleftarrow}}}%
          {\ifdimgreater{\ulen}{\xlen}%
             {\underset{#1}{\overset{#2}{\longleftarrow}}}
             {\xleftarrow[#1]{#2}}}}%
      {\xleftarrow[#1]{#2}}
   }
\newcommand{\isoarrow}{%
   \ifbool{@display}{\overset{\sim}{\longrightarrow}}{\xrightarrow\sim}%
   }
 
\begin{document}

\title[Irreducible components of ADLV in the affine flag variety]{Irreducible components of affine Deligne-Lusztig variety in the affine flag variety}

\author{Arghya Sadhukhan}
\address[A. S.]{Universit\"at M\"unster, Einsteinstr. 62, 48149 M\"unster, Germany}

\email{arghya.s@uni-muenster.de}

\begin{abstract}
We establish a bijection between elements of stable Demazure crystals and labeled longest increasing paths in the double Bruhat graph. As an application, we show that the orbits of irreducible components of certain special affine Deligne–Lusztig varieties (ADLVs) in the affine flag variety of a split reductive group are counted by weight multiplicities of the associated Demazure modules for the dual group. By decomposing an arbitrary ADLV into strata that are iterated fibrations over these special ADLVs, we then obtain a formula for the number of orbits of irreducible components for ADLVs in the affine flag variety, generalizing the numerical form of the Chen–Zhu conjecture.
\end{abstract}

\maketitle

\vspace{-0.8cm}
\section{Introduction}
\subsection{}
Affine Deligne--Lusztig varieties (ADLVs) were introduced by Rapoport
in \cite{Rap05}. In the equal characteristic setting, they are related
to moduli spaces of local shtukas. In the mixed characteristic setting,
they are related to the geometry of Rapoport--Zink spaces and hence,
via $p$-adic uniformization, to the geometry of certain distinguished
loci in the special fiber of Shimura varieties. In this paper, we are interested in the (top-dimensional) irreducible components of
ADLVs: their geometric study
provides useful information about special cycles on
Shimura varieties; we refer the reader to the work on Tate conjectures on Shimura varieties in some special cases in \cite{xiao2017cycles}, \cite{hove2026cycles}.

Let us first recall the corresponding
problem at the affine Grassmannian level. Let $F$ be a local field with the
ring of integers $\mathcal O_F$ and the completion
of its maximal unramified extension $\breve F$. Let $G$ be a reductive group
over $\CO_F$.
For $b\in G(\breve F)$ and a cocharacter $\mu$ of $G$, the affine
Deligne-Lusztig variety $X_\mu(b)$ is a locally closed, equi-dimensional subscheme of
the affine Grassmannian. The scheme $X_\mu(b)$ carries an action of $J_b(F)$, the $F$-rational
points of the reductive group $J_b$ over $F$ given by the
Frobenius centralizer of $b$. This induces an action on
$\Sigma(X_\mu(b))$, the set of irreducible components of $X_\mu(b)$. The resolution of Chen–Zhu conjecture by Nie~\cite{Nie22} and Zhou–Zhu~\cite{ZZ20} then establishes a striking link between the geometry of ADLVs in the affine Grassmannian and the representation theory of its complex Langlands dual group~$\widehat{G}(\mathbb{C})$: the orbit set
$J_b(F)\backslash\Sigma(X_\mu(b))$ is naturally in bijection with the crystal basis of the $\lambda_b$-weight space in the highest weight module~$V(\mu)$ of~$\widehat{G}(\BC)$, where $\lambda_b=\lfloor\nu(b)\rfloor$ encodes the “best integral approximation’’ of the Newton point~$\nu(b)$. 

\subsection{}
In this article, we study the analogous question at the affine flag level for a split reductive group over $F$.
For an element $w$ of the extended affine Weyl group, let $X_w(b)$
denote the corresponding ADLV in the affine
flag variety. Our aim is to understand the $J_b(F)$-set
$\Sigma^{\mathrm{top}}(X_w(b))$ of top-dimensional irreducible components of $X_w(b)$. 

As is well-known, there are many complications that arise in the affine flag case: these ADLVs are not equidimensional in general, nor do we have a complete understanding of their non-emptiness pattern or a dimension formula despite the concentrated efforts of many researchers for the last two decades. However, we do have a fairly good understanding of these ADLVs when the element $w$ is \emph{special}, or is \emph{large} compared to the ``gap" between (translation part of) $w$ and $\nu(b)$. Indeed, under an assumption of the second kind, Schremmer \cite{Sch25} gives an enumeration of top-dimensional irreducible components in terms of certain statistics in the double Bruhat graph (DBG); the latter is a graph whose vertices are elements of the finite Weyl group $W$.
While this does give a complete answer in that restrictive regime, it is a more convoluted combinatorial count than the kind provided by the Chen-Zhu conjecture at the affine Grassmannian level.

\subsubsection{} With the goal of achieving such a simple closed formula in the case of split groups, our first result establishes an interesting connection between certain paths in DBG and a subset of the crystal basis associated to the negative quantum group $\CU_q^-(\hat{\mathfrak{g}})$ of the dual complex Lie algebra: this relies on explicit polyhedral realization of Demazure cystals carried out in Fujita's work \cite{Fu22}.

\begin{theorem}\label{thm:A}
    Fix a reduced word $\mathbf{i}$ for $w_0$ and an associated reflection order $\prec$ on $\Phi^+$. Let $x\in W$; then there is a weight-negating bijection \[\Psi_{\ii}: B_x(\infty) \xlongrightarrow[]{\sim} \mathcal{P}_{\max}^{\prec}(1\rightarrow x).\]

    In particular, we have $m_x^{\text{st}}(\delta)=
 \#\bigl\{(p,m) \in\CP_{\max}^{\prec}(1\rightarrow x): \wt(p,m)=\delta\bigr\},$ for every $\delta\in \BZ_{\geq 0} \Phi^{\vee+}$.
\end{theorem}
Here $\Psi_\ii$ is related to the \emph{Lusztig parametrization} $\Upsilon_\ii$ for the crystal basis $B(\infty)$, restricted to the subset of \emph{stable} Demazure crystal $B_x(\infty)$ associated to a finte Weyl group element $x$, and $\mathcal{P}_{\max}^{\prec}(1\rightarrow x)$ denotes the set of all labeled longest increasing paths in DBG from $1$ to $x$, see Section~\ref{sec:crystalprelim} and Section~\ref{sec:DBGprelim}; both sides have a notion of \emph{weight}, and $m_x^{\text{st}}(\delta)$ is the multiplicity of the $-\delta$ weight space in $B_x(\infty)$. This can be related to the weight space $V_x(\lambda)_{\lambda-\delta}$ in a finite dimensional Demazure module $V_x(\lambda)$ for the dual group as long as $\lambda$ is sufficiently ``large" compared to $\delta$, see Lemma~\ref{lem:demazure-stability}. 
\subsubsection{}
We can now combine Theorem~\ref{thm:A} with Schremmer's results in \cite{Sch25} to deduce information on the irreducible components of ADLVs associated to certain special elements of $\widetilde{W}$, namely of the form $w=xt^\lambda$, where $x\in W$ and $\lambda$ is dominant. It is known that the ADLVs associated with such elements are exceptionally well-behaved: for instance, they are equidimensional with a known dimensional formula, see \cite{VM20}, \cite{He21b}. 

\begin{theorem}\label{thm:B}
    Assume that $G$ is split. Let $w=xt^\lambda \in \widetilde{W}$ where $\lambda$ is dominant and $b\in B(G)$ such that $X_w(b)\neq \emptyset$. Assume further that $\text{depth}(\lambda)\geq \<3\rho,\lambda-\nu(b)\>$.
     Then 
    \[\#~ J_b(F)\backslash \Sigma(X_w(b)) = \dim V_x(\lambda)_{\nu(b)},\]
    the multiplicity of $\nu(b)$ weight space in the Demazure module $V_x(\lambda)$.
\end{theorem}

Here $\text{depth}(\lambda):=\min\limits_{\alpha~\text{simple}} \<\alpha,\lambda\>$ is an estimate of how large $\lambda$ is, and the condition on depth is directly inherited from Schremmer's work. Note that if $x=w_0$, our formula recovers the numerical count from the Chen-Zhu conjecture (albeit under the restrictive depth hypothesis), since $V_{w_0}(\lambda)=V(\lambda)$ and the irreducible components of $X_{w_0t^\lambda}(b)$ are known to be in bijection with those of $X_{\lambda}(b)$ as well, see \cite{HZZ24}[Proposition 2.4.10].

We stress that our aim here is to provide a \emph{representation-theoretic} reformulation of a known answer (in terms of DBG statistics) to the question of irreducible components in the affine flag case in line with their affine Grassmannian counterpart, as well as to predict the shape of an answer beyond the restrictive regime built into the proof via the DBG mechanism. Another purpose is to potentially open up new connections between the study of ADLV and combinatorics of crystal basis. 

\subsubsection{}
Let us now turn to ADLVs $X_w(b)$ for more general elements $w\in\widetilde W$. In \cite{He21b}, He established a
nonemptiness criterion in a broad range of cases by relating $w$ to an element of the form $a t^\gamma$, with $a\in W$ and $\gamma$
dominant. By enhancing this construction and keeping track of the intermediate subschemes of certain products of affine flag varieties arising in this procedure, we relate the top-dimensional irreducible components of $X_w(b)$ to those of the ADLV's associated with such special elements. Let $h$ be the maximum of the Coxeter numbers of the irreducible components of $G$.

\begin{theorem}\label{thm:C-strong}
Assume that $G$ is split. Let
$w=x t^\lambda y\in\widetilde W$, where $x,y\in W, \lambda\in X_*(T)^+,$ with $\operatorname{depth}(\lambda)\geq 2h$. Let
$[b]\in B(G)$ be such that $X_w(b)\neq\emptyset$.
\begin{enumerate}[(i)]
    \item There is an explicit finite family of elements
    $\bigl\{a_{\mathcal P}t^{\gamma_{\mathcal P}}
    :\mathcal P\in\mathscr P(w)_b^{\mathrm{top}}\bigr\},$
    where $a_{\mathcal P}\in W$ and $\gamma_{\mathcal P}$ is dominant,
    and a $J_b(F)$-equivariant bijection
    \[
    \Sigma^{\mathrm{top}}(X_w(b))
    \;\xlongrightarrow{\sim}\;
    \bigsqcup_{\mathcal P\in\mathscr P(w)_b^{\mathrm{top}}}
    \Sigma\bigl(X_{a_{\mathcal P}t^{\gamma_{\mathcal P}}}(b)\bigr).
    \]

    \item If, in addition, $\operatorname{depth}(\lambda)
    \geq 2h+\langle3\rho,\lambda-\nu(b)\rangle,$
    then
    \[
    \#~J_b(F)\backslash
    \Sigma^{\mathrm{top}}(X_w(b))
    =
    \sum_{\mathcal P\in\mathscr P(w)_b^{\mathrm{top}}}
    \dim V_{a_{\mathcal P}}(\gamma_{\mathcal P})_{\nu(b)}.
    \]
\end{enumerate}
\end{theorem}

The index set $\mathscr P(w)_b^{\mathrm{top}}$ records certain reduction
data that contribute to the top-dimensional components of $X_w(b)$. We refer to Section~\ref{subsec:cordial ends} for the
construction of the index set and the associated elements, as well as discussions on when we can weaken the depth hypothesis.

We highlight a situation in which the union has a single summand.
A distinguished candidate for the aforementioned elements is $\eta(w)t^{\gamma_0} \text{ with }
\eta(w)=yx,
\gamma_0=\lambda-\mu+(yx)^{-1}\mu$,
where $\mu$ is a specified cocharacter depending on $y$; when $y=1$,
we have $\mu=0$. In
Proposition~\ref{prop:nested-cover-criterion}, we give an explicit condition under which this endpoint is indexed by the unique
element of $\mathscr P(w)_b^{\mathrm{top}}$. In that case, the bijection
in Theorem~\ref{thm:C-strong} becomes
\[
\Sigma^{\mathrm{top}}(X_w(b))
\;\xlongrightarrow{\sim}\;
\Sigma\bigl(X_{\eta(w)t^{\gamma_0}}(b)\bigr).
\]
Under the additional regularity assumption in
Theorem~\ref{thm:C-strong}(ii), it follows that
\[
\#~J_b(F)\backslash\Sigma^{\mathrm{top}}(X_w(b))
=
\dim V_{\eta(w)}(\gamma_0)_{\nu(b)}.
\]

The proof is related in spirit to Deligne-Lusztig reduction, but uses
cordial elements as endpoints of the comparison, as Theorem~\ref{thm:B} then provides the
representation-theoretic answer to the question at hand.
Moreover, the cardinality of
$\mathscr P(w)_b^{\mathrm{top}}$ is bounded independently of
$\lambda$ in a precise sense.

\subsection{Acknowledgements} This work started while the author was still affiliated with the National University of Singapore and was written up during his time at the Universit\"at M\"unster. He thanks both institutes for their excellent working conditions. He is also grateful to Ian Gleason for encouraging him to think about these questions, and to Marc Besson for helpful suggestions during their collaboration on a related forthcoming project. The author used OpenAI’s GPT-6 Astra to assist only with
proofreading the manuscript and to obtain feedback on possible omissions or errors in the proofs. During this work, the author was partially funded by the Deutsche Forschungsgemeinschaft (DFG, German Research Foundation) under Germany's Excellence Strategy EXC 2044/2–390685587, Mathematics Münster: Dynamics–Geometry–Structure, and byLeibniz prize awarded to Eva Viehmann.

\section{Preliminaries}

Throughout, let $F$ be a non-archimedean local field with residue field
$\mathbb F_q$, and let $\breve F$ be the completion of its maximal
unramified extension. Write $k=\overline{\mathbb F}_q$ for the residue
field of $\breve F$ and $\sigma$ for the Frobenius automorphism of
$\breve F/F$.

Let $G$ be a split connected reductive group over $F$. Fix a split
maximal torus $T\subseteq G$ and a Borel subgroup $B\supseteq T$.
Let $\Phi\subseteq X^*(T)$ be the associated root system, with positive
roots $\Phi^+$ and simple roots
$\Delta=\{\alpha_i:i\in\mathbb S\}$. We write $\alpha_i^\vee\in X_*(T)$ for the simple coroot corresponding to
$\alpha_i$, and put $Q^\vee=\mathbb Z\Phi^\vee,
Q^\vee_+=\sum_{i\in\mathbb S}\mathbb Z_{\geq 0}\alpha_i^\vee,
\rho=\frac12\sum_{\alpha\in\Phi^+}\alpha$.
The set of dominant coweights is
$X_*(T)^+
=
\{\lambda\in X_*(T):
\langle\lambda,\alpha_i\rangle\geq 0
\text{ for all }i\in\mathbb S\}.$
Let $W=N_{\bG}(T)(F)/T(F)$
be the finite Weyl group, with simple reflections
$\{s_i:i\in\mathbb S\}$. We denote its length function by $\ell$ and its
longest element by $w_0$. For $w\in W$, let
$R(w)=
\{(i_1,\ldots,i_{\ell(w)}):
w=s_{i_1}\cdots s_{i_{\ell(w)}}\}$
be the set of reduced words for $w$.

Fix a reduced word $\ii=(i_1,\ldots,i_N)\in R(w_0)$, where
$N=\ell(w_0)$. For $1\leq j\leq N$, set
$\beta_j=s_{i_1}\cdots s_{i_{j-1}}(\alpha_{i_j})\in\Phi^+$.
These roots enumerate $\Phi^+$ and determine the reflection order $\beta_1\prec\beta_2\prec\cdots\prec\beta_N$.

\subsection{Polyhedral realizations of crystal bases}
\label{sec:crystalprelim}

Let $\widehat G$ be the complex Langlands dual group of $G$, with
dual maximal torus $\widehat T$ and a Borel subgroup $\widehat B$
corresponding to $B$. Its based root datum is dual to that of $G$.
In particular, the simple roots and simple coroots of $\widehat G$ are $\widehat\alpha_i=\alpha_i^\vee,
\widehat\alpha_i^\vee=\alpha_i
(i\in\mathbb S)$.
We identify the Weyl group of $\widehat G$ with $W$, so that the simple
reflection associated with $\widehat\alpha_i$ is again denoted by $s_i$.
The positive roots of $\widehat G$ are
$\widehat\Phi^+=(\Phi^+)^\vee$, and the reduced word $\ii$ determines
the reflection order $\beta_1^\vee\prec\beta_2^\vee\prec\cdots\prec\beta_N^\vee$
on $\widehat\Phi^+$. Put $P:=X^*(\widehat T)=X_*(T),
P_+:=X_*(T)^+.$
Thus $P$ is the character lattice of $\widehat T$, and the positive root
lattice of $\widehat G$ is $Q^\vee_+$. Set
$\widehat{\mathfrak g}=\operatorname{Lie}(\widehat G)$.
The crystal constructions below use the semisimple derived Lie algebra
of $\widehat{\mathfrak g}$, with weights lying in $P$.

We use the standard conventions for Kashiwara crystals
\cite{K91}, \cite{K93}. Write $B(\infty)$ for the crystal of the negative half $U_q^- (\widehat{\mathfrak{g}})$ of
the quantized enveloping algebra of $\widehat{\mathfrak g}$, with
distinguished element $b_\infty\in B(\infty)$ of weight zero. For $\lambda\in P_+$,
let $V(\lambda)$ be the irreducible representation of $\widehat G$ of
highest weight $\lambda$, and let $B(\lambda)$ be its crystal, with
highest element $b_\lambda$. The canonical projection $\pi_\lambda:B(\infty)\longrightarrow B(\lambda)\sqcup\{0\}$
restricts to a bijection
$\widetilde B(\lambda):=
\{b\in B(\infty):\pi_\lambda(b)\neq0\}
\xrightarrow{\sim} B(\lambda)$.

We denote the Kashiwara operators by
$\ee_i$ and $\ff_i$. For a reduced expression $w=s_{j_1}\cdots s_{j_r}$, define
\[
\begin{aligned}
B_w(\infty)
&=
\{\ff_{j_1}^{a_1}\cdots\ff_{j_r}^{a_r}b_\infty:
a_1,\ldots,a_r\geq0\},\\
B_w(\lambda)
&=
\{\ff_{j_1}^{a_1}\cdots\ff_{j_r}^{a_r}b_\lambda:
a_1,\ldots,a_r\geq0\}\setminus\{0\}.
\end{aligned}
\]
These sets are independent of the chosen reduced expression, and
$\pi_\lambda$ identifies
$B_w(\infty)\cap\widetilde B(\lambda)$ with $B_w(\lambda)$.

Let $\widehat{\mathfrak b}$ and $\widehat{\mathfrak b}^{-}$ be the
positive and negative Borel subalgebras of $\widehat{\mathfrak g}$.
For $w\in W$, choose a nonzero extremal vector
$v_{w\lambda}\in V(\lambda)_{w\lambda}$. The \emph{Demazure module}
\[
V_w(\lambda)
:=
U(\widehat{\mathfrak b})v_{w\lambda} \subset V(\lambda)
\]
is spanned by the lower global-basis vectors indexed by $B_w(\lambda)$.
Consequently,
\begin{equation}\label{eq:character}
\operatorname{ch}V_w(\lambda)
=
\sum_{b\in B_w(\lambda)}e^{\operatorname{wt}(b)},
\qquad
\dim V_w(\lambda)_\mu
=
\#\{b\in B_w(\lambda):\operatorname{wt}(b)=\mu\}.
\end{equation}
Similarly, the opposite Demazure module $V^v(\lambda):=
U(\widehat{\mathfrak b}^{-})v_{v\lambda}$
is indexed by the opposite Demazure crystal
$B^v(\lambda)\subseteq B(\lambda)$.

For $\delta\in Q^\vee_+$, define the \emph{stable} Demazure multiplicity
\[
m_w^{\mathrm{st}}(\delta)
:=
\#B_w(\infty)_{-\delta}
=
\#\{b\in B_w(\infty):\operatorname{wt}(b)=-\delta\}.
\]
This set is finite. For fixed $w$ and $\delta$, we have $m_w^{\mathrm{st}}(\delta)
=
\dim V_w(\lambda)_{\lambda-\delta}$
whenever $\lambda$ is sufficiently dominant; see Lemma~\ref{lem:demazure-stability}.

Define the involution $i\mapsto i^*$ of $\mathbb S$ by
$w_0(\alpha_i)=-\alpha_{i^*}$, or equivalently by
$w_0(\alpha_i^\vee)=-\alpha_{i^*}^\vee$, and set
$\ii^*=(i_1^*,\ldots,i_N^*)\in R(w_0)$. The Sch\"utzenberger involution $\eta_\lambda:B(\lambda)\to B(\lambda)$
satisfies $\eta_\lambda\ff_i=\ee_{i^*}\eta_\lambda,
\operatorname{wt}(\eta_\lambda(b))
=
w_0\operatorname{wt}(b)$.
It interchanges ordinary and opposite Demazure crystals:
\[
\eta_\lambda(B_w(\lambda))=B^{w_0w}(\lambda),
\qquad
\eta_\lambda(B^v(\lambda))=B_{w_0v}(\lambda).
\]

We recall two parametrizations associated with $\ii$. Following \cite{Lit98}[\S 1], for
$b\in B(\infty)$ or $B(\lambda)$, set $b^{(0)}=b$ and recursively
define
\[
a_j=\max\{m\geq0:\ee_{i_j}^{\,m}b^{(j-1)}\neq0\},
\qquad
b^{(j)}=\ee_{i_j}^{\,a_j}b^{(j-1)}
\quad(1\leq j\leq N).
\]
The resulting injective map
$\Phi_\ii:b\longmapsto(a_1,\ldots,a_N)$
is called the \emph{Berenstein–Littelmann–Zelevinsky string parametrization}. There is a rational polyhedral cone
$C_\ii\subseteq\mathbb R_{\geq0}^N$, called the \emph{string cone}, such that
$\Phi_\ii(B(\infty))=C_\ii\cap\mathbb Z^N$.
For $\lambda\in P_+$, the \emph{string polytope} is
\begin{equation}\label{eq:string}
    \Delta_\ii(\lambda)
=
\left\{
a\in C_\ii:
a_j+\sum_{k>j}
\langle\alpha_{i_k}^\vee,\alpha_{i_j}\rangle a_k
\leq
\langle\lambda,\alpha_{i_j}\rangle
\quad(1\leq j\leq N)
\right\},
\end{equation}
and $\Phi_\ii(B(\lambda))
=
\Delta_\ii(\lambda)\cap\mathbb Z^N$, see \cite{BZ01}[\S 3.2], \cite{Lit98}[\S 1].

The PBW basis associated with $\ii$ gives in \cite{Lus96} \emph{Lusztig's parametrization} $\Upsilon_\ii:B(\infty)\xrightarrow{\sim}\mathbb Z_{\geq0}^N$.
Using the unique lift in $\widetilde B(\lambda)$, we also regard
$\Upsilon_\ii$ as a map on $B(\lambda)$. Its image is the set of lattice
points of a rational convex polytope $\widehat\Delta_\ii(\lambda)$,
called the \emph{Lusztig polytope}:
\begin{equation}\label{eq:lusztigpoly}
\Upsilon_\ii(B(\lambda))
=
\widehat\Delta_\ii(\lambda)\cap\mathbb Z^N.
\end{equation}

The two parametrizations are related by an affine unimodular
transformation. Define
$\Omega_{\ii,\lambda}:\mathbb R^N\to\mathbb R^N$ by
$\Omega_{\ii,\lambda}(a_1,\ldots,a_N)=(a_1',\ldots,a_N')$, where
\[
a_j'
=
\langle\lambda,\alpha_{i_j}\rangle-a_j
-
\sum_{k>j}
\langle\alpha_{i_k}^\vee,\alpha_{i_j}\rangle a_k.
\]
By Morier-Genoud's formula \cite{MG08}[Corollary 2.17],
\begin{equation}\label{eq:MG}
\Omega_{\ii,\lambda}\circ\Phi_\ii
=
\Upsilon_{\ii^*}\circ\eta_\lambda
\quad\text{on }B(\lambda),
\qquad
\Omega_{\ii,\lambda}(\Delta_\ii(\lambda))
=
\widehat\Delta_{\ii^*}(\lambda).
\end{equation}
Thus $\Omega_{\ii,\lambda}$ takes the string coordinates of $b$ with
respect to $\ii$ to the Lusztig coordinates of $\eta_\lambda(b)$ with
respect to $\ii^*$.

\subsection{Double Bruhat graph}
\label{sec:DBGprelim}

We use the double Bruhat graph associated with the root system of
$\bG$, with the reflection order $\prec$ fixed above. This was introduced by Naito-Watanabe \cite{NW17} to describe periodic R-polynomials of affine Weyl groups.

\begin{definition}
The double Bruhat graph $\operatorname{DBG}(W)$ has vertex set $W$.
For each $u\in W$ and $\alpha\in\Phi^+$, it has a directed edge
\[
u\xrightarrow{\ \alpha\ }us_\alpha.
\]

An unlabeled path $p$ from $v$ to $w$ is a sequence of adjacent edges
\[
p:\quad
v=u_0\xrightarrow{\alpha_1}u_1
\xrightarrow{\alpha_2}\cdots
\xrightarrow{\alpha_r}u_r=w.
\]
It is \emph{increasing} with respect to $\prec$ if
$\alpha_1\prec\cdots\prec\alpha_r$.

A labeled path $(p,m)$ consists of such a path and integers
$m=(m_1,\ldots,m_r)$ satisfying
\[
m_j\geq
\begin{cases}
0,&\ell(u_j)>\ell(u_{j-1}),\\
1,&\ell(u_j)<\ell(u_{j-1}).
\end{cases}
\]
Its length and weight are
\[
\ell(p)=\ell(p,m):=r,
\qquad
\operatorname{wt}(p,m):=
\sum_{j=1}^r m_j\alpha_j^\vee\in Q^\vee_+.
\]
It is increasing if its underlying unlabelled path is increasing.

For $x,y \in W$, we write $\CP^\prec_{\max}(x\to y)$, respectively
$\bar{\CP}^\prec_{\max}(x\to y)$, for the set of longest increasing
labeled paths, respectively their underlying unlabeled paths, from $x$
to $y$. We omit decorations when they are clear from the context.
\end{definition}

By \cite{Sch25}[Lemma~4.11],
\[
\max\{\ell(p):p\text{ is an increasing path from }x\text{ to }y\}
=
\ell(x^{-1}y).
\]
In particular, every path in $\CP^\prec_{\max}(1\to x)$ has length
$\ell(x)$.

Under the dual-root-datum identification from Section~\ref{sec:crystalprelim}, the weight of a labeled path is an element of
the positive root lattice of $\widehat G$. Thus the path weight
$\operatorname{wt}(p,m)$ and the crystal weight
$-\operatorname{wt}(b)$ for $b\in B(\infty)$ lie in the same lattice
$Q^\vee_+$.

\subsection{Affine Deligne--Lusztig varieties}

Let $\breve G=G(\breve F)$. We choose the base alcove $\mathfrak a$ in the dominant chamber and
normalize the action on the apartment so that $t^\lambda=\lambda(t)$ acts
by translation by $\lambda$. Thus the $\s$-stable Iwahori subgroup $I$ fixing
$\mathfrak a$ is the inverse image of the opposite Borel $B^-$ under
$G(\mathcal O)\to G(k)$. The corresponding positive affine roots,
viewed as affine functions on the apartment, are
\[
\Phi_{\mathrm{aff}}^+
=\Phi^+\sqcup
\{\alpha+n\delta_{\mathrm{aff}}:\alpha\in\Phi,\ n\geq1\};
\]
in particular, if $\Phi$ is irreducible, the simple affine roots are
$\Delta\cup\{-\theta+\delta_{\mathrm{aff}}\}$, where $\theta$ is the
highest root. Choosing a special vertex in its
closure identifies the Iwahori-Weyl group with
$\tW
=
N_{G}(T)(\breve F)/
\bigl(T(\breve F)\cap\breve\CI\bigr)
\cong X_*(T)\rtimes W$. We may thus express an arbitrary element $w$ as $xt^\lambda y$, where $x,y\in W$, $\lambda$ is dominant and $t^\lambda y$ is a minimal length element in $W\backslash \tW$; with the above conventions, we then have
\[
\ell(xt^\lambda y)
=\ell(x)+\langle2\rho,\lambda\rangle-\ell(y).
\]
With such a representation of $w$, we also define $\eta(w):=yx \in W$. The $\sigma$-conjugation action of $\breve G$ on itself is given by $g\cdot_\sigma b=gb\sigma(g)^{-1}$. The Frobenius $\sigma$ acts trivially on $X_*(T)$, $W$, and $\tW$. Write $[b]$ for the $\sigma$-conjugacy class of $b\in\breve G$ and
$B(G)$ for the set of such classes. By Kottwitz's classification in \cite{Ko85} and \cite{Ko97},
$[b]$ is determined by its \emph{Kottwitz point} and dominant \emph{Newton point}
\[
\kappa(b)\in\pi_1(G):=X_*(T)/Q^\vee,~
\nu(b)\in X_*(T)_{\mathbb Q}^+.
\]
For $\lambda\in X_*(T)$, let $\lambda^\natural$ denote its image in
$\pi_1(G)$. If $\lambda$ is dominant, then
$\kappa_G(t^\lambda)=\lambda^\natural,
\nu(t^\lambda)=\lambda$.

For $w\in\tW$, choose a representative
$\dot w\in N_{G}(T)(\breve F)$. The \emph{affine Deligne-Lusztig variety}
associated with $w$ and $b$ is
\[
X_w(b)
=
\{g\CI\in\breve G/\CI:
g^{-1}b\sigma(g)\in\CI\dot w\CI\}.
\]
In equal characteristic, we regard $X_w(b)$ as a reduced locally closed
subscheme of the affine flag variety. In mixed characteristic, we regard
it as a locally closed perfect subscheme of the $p$-adic affine flag
variety; see \cite{Zhu17}, \cite{BS17}.

For $w\in \widetilde{W}$, we define $B(G)_w:=\{[b]: X_w(b)\neq \emptyset\} \subset B(G)$. Let $J_b$ be the $\sigma$-centralizer of $b$, so that $J_b(F)
=
\{g\in\breve G:g^{-1}b\sigma(g)=b\}$. We define the defect of $b$ as $\operatorname{def}_G(b)
=
\operatorname{rank}_F(G)-\operatorname{rank}_F(J_b).$ Following He \cite[\S 10.1]{He14}, we define virtual dimension for the pair $(w,b)$ as
\[d_w(b):= \frac{1}{2} \{\ell(w)+\ell(\eta(w))-\<2\rho,\nu(b)\>-\mathrm{def}(b)\}.\]
By corollary 10.4 in loc. cit., we then have $\dimm X_w(b)\leq d_w(b)$ whenever $\kappa(w)=\kappa(b)$. 

Let us also note, if $G= G_1 \times G_2$ then $B(G)=B(G_1)\times B(G_2)$
and $\widetilde{W}_G=\widetilde{W}_{G_1}\times\widetilde{W}_{G_2}$ . Hence, $X_w(b)\simeq X_{w_1}(b_1)\times X_{w_2}(b_2)$ where $w_i$’s are
the projections of $w$ onto $\widetilde{W}_{G_i}$ and $b_i$’s are that of $b$ onto $B(G_i)$. Moreover, $J_b(F)=J_{b_1}(F)\times J_{b_2}(F)$, and the
top-dimensional irreducible components of $X_w(b)$ are the
products of those of the two factors. On the representation-theoretic
side, the dual group decomposes as
$\widehat G=\widehat G_1\times\widehat G_2$, and the corresponding
Demazure modules decompose as $V_{(a_1,a_2)}(\lambda_1,\lambda_2)
\simeq
V_{a_1}(\lambda_1)\boxtimes V_{a_2}(\lambda_2)$.
In particular, their weight multiplicities factor as products. Thus our component-counting results for a direct product follow
from the corresponding results for its factors.

The group $J_b(F)$ acts on $X_w(b)$ by left multiplication. We write
$\Sigma^{\mathrm{top}}(X_w(b))$ for its set of top-dimensional
irreducible components; for equidimensional ADLVs, we may drop the decoration ``top". Although this set may be infinite, the orbit set $J_b(F)\backslash\Sigma^{\mathrm{top}}(X_w(b))$ is known to be finite by \cite{HV20}. 

For $\lambda\in X_*(T)^+$, define
$\text{depth}(\lambda)
=
\min_{i\in\mathbb S}\langle\lambda,\alpha_i\rangle$. For $x\in W$, let
$\supp(x)\subseteq\mathbb S$ be the set of simple
reflections appearing in a(ny) reduced expression for $x$. We will be particularly interested in elements of the form
$w=xt^\lambda$, where $x\in W$ and $\lambda$ is dominant. By
\cite[Theorem~4.2]{He21b}, such elements belong to the family of so-called \emph{cordial elements} \`a la Milićević-Viehmann \cite{VM20}. Combining this with the
results of Milićević-Viehmann on cordial elements, Viehmann's
description of the minimal Newton point, and the comparison with Levi
subgroups gives the following description for the non-emptiness and dimension in the regular range, specialized to the case of a split group; see
\cite[Theorem~1.1]{VM20}, \cite[Theorem~1.1(1)]{Vie21}, and
\cite[Theorem~2]{Sch24}.

\begin{theorem}\label{prop:xtlambda}
Assume that $G$ is split. Let $w=xt^\lambda$ where $\lambda\in X_*(T)^+, x\in W$ and $J=\supp(x)$.
\begin{enumerate}[(i)]
\item  Assume that $\lambda$ satisfy
$\operatorname{depth}(\lambda)\geq2$. Then
\[
X_w(b)\neq\varnothing
\quad\Longleftrightarrow\quad
\kappa_G(b)=\lambda^\natural
\quad\text{and}\quad
\lambda-\nu(b)\in
\sum_{i\in J}\mathbb Q_{\geq0}\alpha_i^\vee.
\]
\item Whenever $X_w(b)$ is nonempty, it is equidimensional of dimension
\[
\dim X_w(b)
=
\langle\lambda-\nu(b),\rho\rangle
+\ell(x)-\frac12\operatorname{def}_G(b).
\]
\end{enumerate}
\end{theorem}

Using the language of DBG, Schremmer gives the following enumeration for the irreducible components for the ADLVs associated to such cordial elements in \cite[Corollary 5.9]{Sch25}.

\begin{theorem}\label{thm:DBGcount}
Assume that $G$ is split. Let $w=xt^\lambda$ where $\lambda\in X_*(T)^+, x\in W$. Suppose that
$X_w(b)\neq\emptyset$ and that $\text{depth}(\lambda) \geq 3\langle\lambda-\nu(b),\rho\rangle$.
Then
\[
\#~
J_b(F)\backslash
\Sigma^{\mathrm{top}}(X_w(b))
=
\#\left\{
(p,m)\in\CP^\prec_{\max}(1\rightarrow x):
\operatorname{wt}(p,m)=\lambda-\nu(b)
\right\}.
\]
\end{theorem}

We note that the stated depth assumption implies in particular that
$\nu(b)$ is regular. Since $G$ is split, $\nu(b)$ is then integral and
$\operatorname{def}_G(b)=0$. In particular,
$\lambda-\nu(b)\in Q^\vee_+$.
\section{Double Bruhat graph and Demazure module}
\subsection{Union of faces}
We first recall the explicit description of Berenstein–Littelmann–Zelevinsky string polytope and Lusztig polytope of a Demazure module. 

For $\kk\subset[N]$, define
\[
F_\kk(\Delta_\ii(\lambda))=\{a\in\Delta_\ii(\lambda):
a_j+\sum_{k>j} \<\alpha_{i_k}^\vee,\alpha_{i_j}\>a_k=\<\lambda,\alpha_{i_j}\>~ \forall j\in\kk\}.
\]

Similarly, set $E_\kk=\{x\in\R_{\ge0}^N:x_j=0\ ~ \forall j\in\kk\}$. Following Fujita \cite[\S 4]{Fu22}, we define, for $v\in W$, the set
\[R(\ii,v) := \{\kk = (k_1,\cdots,k_{\ell(v)}) : 1 \leq k_1 < \cdots <k_{\ell(v)} \leq N, (i_{k_1},\cdots,i_{k_{\ell(v)}}) \in R(v)\}.\]

As noted in loc. cit., they can be viewed as a generalization of reduced pipe dreams for arbitrary Lie type. We set $R(\ii,e)=\{\emptyset\}$.  

\begin{theorem}\cite[Theorem 4.5 and Corollary 4.6]{Fu22}\label{thm:fujita}
For every $v,w\in W$ and every dominant integral $\lambda$, one has
\begin{align*}
\Phi_\ii(B^v(\lambda))
&=\left(\bigcup_{\kk\in R(\ii,v)}F_\kk(\Delta_\ii(\lambda))\right)\cap\Z^N,
\\
\Upsilon_\ii(B_w(\lambda))
&=\left(\widehat\Delta_\ii(\lambda)\cap \bigcup_{\kk\in R(\ii^*,w_0w)}E_\kk \right) \cap\Z^N.
\end{align*}
Moreover, the corresponding statement for the infinite Demazure
crystal is
\[
\Upsilon_\ii(B_w(\infty))
=
\left(
\bigcup_{\kk\in R(\ii^*,w_0w)}E_\kk
\right)\cap\Z^N.
\]
\end{theorem}

\begin{proof}
The first equality is \cite[Theorem 1]{Fu22}. Note that we obtain directly from \eqref{eq:MG},
\begin{equation}\label{eq:faces}
\Omega_{\ii,\lambda}(F_\kk(\Delta_\ii(\lambda)))
=\widehat\Delta_{\ii^*}(\lambda)\cap E_\kk.
\end{equation}
Applying the first equality now with $\ii^*$ and $v=w_0w$, and then using
\eqref{eq:MG}, \eqref{eq:faces}, and
$\eta_\lambda(B^{w_0w}(\lambda))=B_w(\lambda)$, we obtain
\begin{align*}
\Upsilon_\ii(B_w(\lambda))
&=\Omega_{\ii^*,\lambda}
  \bigl(\Phi_{\ii^*}(B^{w_0w}(\lambda))\bigr)\\
&=\left(
  \widehat\Delta_\ii(\lambda)\cap
  \bigcup_{\kk\in R(\ii^*,w_0w)}E_\kk
  \right)\cap\mathbb Z^N.
\end{align*}

For the final assertion, identify $B(\lambda)$ and
$B_w(\lambda)$ with their canonical lifts into $B(\infty)$.
By \cite[Corollary~4.4.5]{K91}, these exhaust
$B(\infty)$ and $B_w(\infty)$, respectively, as $\lambda$
ranges over $P^+$. In particular,
\[
\bigcup_{\lambda\in P^+}
\bigl(\widehat\Delta_\ii(\lambda)\cap\mathbb Z^N\bigr)
=\Upsilon_\ii(B(\infty))
=\mathbb Z_{\ge0}^N.
\]
Taking the union over $\lambda$ in the second equality
therefore proves the final assertion.
\end{proof}

Let $x\in W$. For $p \in \bar{\CP}^{\prec}_{\max}(1 \rightarrow x)$, we define $I(p)
:=\{k_1< \cdots < k_{\ell(x)}\}\subset [N]$ such that \begin{align*}
p : v_0=1\xrightarrow{\beta_{k_1}} v_1\xrightarrow{\beta_{k_2}} \cdots\xrightarrow{\beta_{k_{\ell(x)}}} v_{\ell(x)}=x.
\end{align*}
In view of Theorem~\ref{thm:fujita}, the connection between Demazure crystals and the double Bruhat graph stems from the following simple observation.
\begin{lemma}\label{lem:complement}
We have a bijection 
 \[\mathrm{Del}: \bar{\CP}^{\prec}_{\max}(1 \rightarrow x) \xrightarrow{} R(\ii, xw_0),~ p \mapsto [N]\setminus I(p). \]
\end{lemma}

\begin{proof}

To reduce notational burden, we write $t_k=s_{i_k}$ where $w_0=s_{i_1}\cdots s_{i_N}$ is our fixed word. Let $P\subseteq[N]$ and $K=[N]\setminus P$. Define

\[r_j=\prod_{k\leq j,\,k\in P}s_{\beta_k}; u_j = s_{i_1}\cdots s_{i_j};
q_j=\prod_{k\leq j,\,k\in K}t_k, \text{ with } r_0=q_0=1.\]

We first prove:

(a) $r_ju_j=q_j$ for every $j$; in particular, $\prod_{j\in P}s_{\beta_j}
 =\left(\prod_{k\in K}t_k\right)w_0.$
 
Induct on $j$. If $j\in K$, then
$r_ju_j=r_{j-1}u_{j-1}t_j=q_{j-1}t_j=q_j$.
If $j\in P$, use $s_{\beta_j}=u_{j-1}t_ju_{j-1}^{-1}$ to obtain
$r_ju_j=r_{j-1}u_{j-1}t_j^2=q_{j-1}=q_j$.
Taking $j=N$ proves the final claim in (a).

Now letting $P=I(p)$, we get 

(b) $x=\left(\prod_{k\in K}t_k\right)w_0$, i.e. $xw_0=\prod_{k\in K}t_k$. 

Note that every increasing path $p$ in DBG starting at $1$ is uniquely specified by the set $P:=I(p)$
of root positions; by (b), it ends at $x$ exactly when
the \emph{complementary word}, i.e. the one obtained from $\ii$ by deleting the simple reflections at positions in $I(p)$, has product $xw_0$. Consequently
$N-|P|=|K|\geq\ell(xw_0)=N-\ell(x)$, so $|P|\leq\ell(x)$; equality $\ell(p)=|P|=\ell(x)$ holds precisely when the
complementary word is reduced.
\end{proof}
\begin{lemma}
For every $v\in W, R(\ii,v)=R(\ii^*,w_0vw_0)$. In particular, we have $R(\ii,xw_0)=R(\ii^*,w_0x)$.
\end{lemma}

\begin{proof}
For any increasing tuple $\kk=(k_1,\ldots,k_r)$, we have
\[
s_{i_{k_1}^*}\cdots s_{i_{k_r}^*}
=w_0\bigl(s_{i_{k_1}}\cdots s_{i_{k_r}}\bigr)w_0.
\]
Since conjugation by $w_0$ permutes the simple reflections,
it preserves length and reducedness. Hence
$\kk\in R(\ii,v)$ if and only if
$\kk\in R(\ii^*,w_0vw_0)$.
The final assertion follows by taking $v=xw_0$.
\end{proof}

While Lemma~\ref{lem:complement} provides us with a reinterpretation of the index set in terms of \emph{unlabeled} longest paths in DBG, we note that this union above is not disjoint: this is a problem insofar as we want to relate the \emph{labeled} longest paths - where the labels are going to be supplied by the summand corresponding to the paths - to Demazure crystals. In the next subsection we resolve this issue by choosing a finer decomposition, which also shows why the admissibility condition on labels in DBG paths is natural.

\subsection{Disjointness of union} As mentioned already, a point of the union can belong to several coordinate faces: this happens when its zero set contains more than one reduced subword for $xw_0$. To obtain a disjoint decomposition, we must choose one of these subwords \emph{canonically}. The following lemma makes this choice by selecting the rightmost possible reduced subword, in the sense made precise below.
\begin{lemma}\label{lem:greedy}
Let $v\in W$ and assume that $Z\subseteq[N]$ contains an element of $R(\ii,v)$ as a subset. There is a unique $K\subseteq Z$ in $R_\ii(v)$ such that, writing
$q_j=\prod_{k\leq j,\,k\in K}t_k$, one has
\begin{equation}\label{eq:omitted}
 \ell(q_{j-1}t_j)>\ell(q_{j-1})
 \qquad\text{for every }j\in Z\setminus K.
\end{equation}
It is recovered by scanning positions from $N$ down to $1$, starting with
$q_N=v$, and selecting $j$ exactly when $j\in Z$ and
$\ell(q_jt_j)<\ell(q_j)$.
\end{lemma}

\begin{proof}
For existence, choose $K\subseteq Z$ in $R_\ii(v)$ by maximizing
$\Xi_K:=\sum_{k\in K}2^k$. Suppose that an omitted position $j\in Z\setminus K$
violates \eqref{eq:omitted}. Write $q_{j-1}$ for the reduced subword of $v$ obtained before this stage, i.e.
$q_{j-1}=t_{k_1}\cdots t_{k_a}$, where $k_1<\cdots<k_a<j$.
The exchange condition gives an index $h$ with
$q_{j-1}t_j=t_{k_1}\cdots\widehat{t_{k_h}}\cdots t_{k_a}$, and thus $q_{j-1}=t_{k_1}\cdots\widehat{t_{k_h}}\cdots t_{k_a}t_j$.
Multiplying this equality on the right by the remaining simple reflections occurring in $K$ in order  shows that replacing
the selected position $k_h$ by $j$ leaves the full selected product equal
to $v$, and the resulting word still has $\ell(v)$ letters so it is reduced.
However, $\Xi_{K\cup\{j\}\setminus \{k_h\}}>\Xi_K$, a contradiction.

For uniqueness, consider any $K$ satisfying the conclusion. If $j\in K$,
then its selected prefixes are reduced, so $q_j=q_{j-1}t_j$ and
$\ell(q_jt_j)<\ell(q_j)$. If $j\in Z\setminus K$, then $q_j=q_{j-1}$
and \eqref{eq:omitted} gives $\ell(q_jt_j)>\ell(q_j)$.
Positions outside $Z$ cannot be selected. Thus $q_N=v$ determines every
preceding state and every selected position by the forced recursion
\begin{equation}\label{eq:greedy}
 q_{j-1}=
 \begin{cases}
   q_jt_j,&j\in Z\text{ and }\ell(q_jt_j)<\ell(q_j),\\
   q_j,&\text{otherwise}.
 \end{cases}
\end{equation}
This proves uniqueness and the algorithm.
\end{proof}

\begin{remark}
If $Z\in R_{\ii}(v)$, then any $K\subseteq Z$ belonging to $R_{\ii}(v)$ must equal $Z$, since both have cardinality $\ell(v)$. In this case \eqref{eq:omitted} is  vacuous. The content of Lemma~\ref{lem:greedy} concerns a set $Z$ which merely \emph{contains} a reduced subword for $v$ and may contain several such subwords: the omitted-position condition selects a unique one. For example, in type $A_2$, take $\mathbf i=(1,2,1)$, $v=s_1$, and $Z=\{1,3\}$. The two candidates are $K_1=\{1\}$ and $K_2=\{3\}$. For $K_1$, position $3$ is omitted and $q_2=s_1$, so $\ell(q_2t_3)=\ell(s_1s_1)=0<1=\ell(q_2)$, violating \eqref{eq:omitted}. For $K_2$, the only omitted position in $Z$ is $1$, and $\ell(q_0t_1)=\ell(s_1)=1>0=\ell(q_0)$. Thus $K_2=\{3\}$ is the unique permitted choice, exactly as prescribed by the backwards algorithm. 

As we shall see below in our application, $Z=Z(y):=\{j:y_j=0\}$ is the zero-coordinate set of a point $y$ in the union of certain coordinate faces. Membership in this union means that $Z(y)$ contains some $K\in R_{\ii}(v)$, not necessarily that $Z(y)$ itself belongs to $R_{\ii}(v)$. If $Z(y)=K$, the point lies in the relative interior of the corresponding maximal face and its assignment is already unique. When additional coordinates vanish, the point may lie in several maximal faces, and the lemma resolves this ambiguity.
\end{remark}
For $p\in\bar{\CP}$, we define the set of downward-edge positions by
\[
D(p):=
\bigl\{
k_a\in I(p):
\ell(v_a)<\ell(v_{a-1})
\bigr\}.
\]
Thus $D(p)$ records the positions in the reflection order of the
roots labeling the downward edges of $p$. Define the partially open coordinate orthant
\[
H_p:=
\bigl\{
y\in E_{\mathrm{Del}(p)}:
y_j>0\text{ for every }j\in D(p)
\bigr\}.
\]

\begin{lemma}\label{lem:partition}
With the preceding notation,
\begin{equation}\label{eq:partition}
\bigcup_{p \in \bar{\CP}} E_{\mathrm{Del}(p)}
=\bigsqcup_{p\in\bar{\CP}}H_p.
\end{equation}
\end{lemma}

\begin{proof}
Fix $p\in\bar{\CP}$, put $P=I(p)$ and $K=[N]\setminus P$, and use the
prefix notation of Lemma~\ref{lem:complement}. For $j\in P$, the path
vertex immediately before the edge labeled $\beta_j$ is $r_{j-1}$.
Since $r_{j-1}\beta_j
=
r_{j-1}u_{j-1}\alpha_{i_j}
=
q_{j-1}\alpha_{i_j}$,
the edge at position $j$ is downward precisely when
\begin{equation}\label{eq:down}
j\in D(p)
\quad\Longleftrightarrow\quad
q_{j-1}\alpha_{i_j}\in-\Phi^+
\quad\Longleftrightarrow\quad
\ell(q_{j-1}t_j)<\ell(q_{j-1}).
\end{equation}

Take $y\in 
\bigcup_{\kk\in R(\ii,xw_0)}E_\kk$ and let $Z=\{j:y_j=0\}$. By definition,
$Z$ contains a reduced subword for $xw_0$. Apply Lemma~\ref{lem:greedy}
with $v=xw_0$ to obtain $K\subseteq Z$, and let $P=[N]\setminus K$.
By Lemma~\ref{lem:complement} we then associate to this $P$ an unique path $p\in\overline \CP$, and
$y\in E_K=E_{\text{Del}(p)}$. If $j\in P$ has $y_j=0$, then $j\in Z\setminus K$;
\eqref{eq:omitted} and \eqref{eq:down} show that this edge is upward.
Hence every downward coordinate is positive, so $y\in H_p$.

Conversely, suppose $y\in H_p$. The complement $K=[N]\setminus I(p)$ is a reduced subword
for $xw_0$ contained in $Z$. Each $j\in Z\setminus K$ is a path position with \(y_j=0\), hence is upward by the definition of $H_p$. By
\eqref{eq:down}, $K$ satisfies \eqref{eq:omitted}. Lemma~\ref{lem:greedy}
therefore determines $K$, and consequently $p$, uniquely from $y$.
This proves both coverage and disjointness.
\end{proof}
\begin{example}
We give an example of the above decomposition in type $A_3$. Take $x=s_2s_3s_2s_1$ and $\ii=(1,2,3,1,2,1)$.
Then the ordered roots are

\[\alpha_1 \prec \alpha_1+\alpha_2 \prec \alpha_1+\alpha_2+\alpha_3 \prec
\alpha_2 \prec \alpha_2+\alpha_3 \prec \alpha_3.\]

The reduced subword position sets for $xw_0=s_1s_2$ are
$K_A=\{4,5\}$, $K_B=\{1,5\}$ and $K_C=\{1,2\}$.
Their complements are $P_A=\{1,2,3,6\}$, $P_B=\{2,3,4,6\}$ and
$P_C=\{3,4,5,6\}$, respectively. Using \eqref{eq:down}, the downward
position sets are $D_A=\varnothing$, $D_B=\{4\}$ and $D_C=\{5\}$.
Consequently, the decomposition is
\begin{equation}\label{eq:a3partition}
 E_{K_A}\cup E_{K_B}\cup E_{K_C} = E_{K_A}\;\sqcup\;
       (E_{K_B}\cap\{y_4>0\})\;\sqcup\;
       (E_{K_C}\cap\{y_5>0\}).
\end{equation}
The Lusztig cone $\Upsilon_\ii(B_x(\infty))$ consists of the points with integer coordinates in either of the sides. Now, the increasing unlabeled paths of maximal length \(4\) from \(1\) to \(x\) involve the following ordered tuples of positive roots
\[
p_A:
\quad
\alpha_1,\ 
\alpha_1+\alpha_2,\ 
\alpha_1+\alpha_2+\alpha_3,\ 
\alpha_3,
\]
\[
p_B:
\quad
\alpha_1+\alpha_2,\ 
\alpha_1+\alpha_2+\alpha_3,\ 
\alpha_2,\ 
\alpha_3,
\]
\[
p_C:
\quad
\alpha_1+\alpha_2+\alpha_3,\ 
\alpha_2,\ 
\alpha_2+\alpha_3,\ 
\alpha_3.
\]

The successive Weyl-group lengths along these paths are
\[
\begin{array}{c|c}
\text{path} & \text{successive lengths} \\ \hline
p_A & 0,1,2,3,4;\\
p_B & 0,3,4,3,4;\\
p_C & 0,5,6,3,4.
\end{array}
\]
Thus every edge of \(p_A\) is length-increasing, whereas the third edge of each of \(p_B\) and \(p_C\) is length-decreasing. Accordingly, the label conditions in DBG are
\begin{equation}\label{eq:admissible-labels}
\begin{aligned}
p_A:\quad & m_1,m_2,m_3,m_4\geq 0;\\
p_B:\quad & m_1,m_2,m_4\geq 0,\qquad m_3\geq 1;\\
p_C:\quad & m_1,m_2,m_4\geq 0,\qquad m_3\geq 1.
\end{aligned}
\end{equation}
Noting that the positive coordinate $y_4$ is the $\alpha_2$-label in $p_B$, and the positive coordinate $y_5$ is the $\alpha_2+\alpha_3$-label in $p_C$, we conclude that the disjoint decomposition in Eq.~\ref{eq:a3partition} exactly captures the admissibility of labels in DBG in Eq.~\ref{eq:admissible-labels}.
\end{example}
\subsection{Truncating to finite levels}
We start by recording the weight of crystal elements in both parametrizations.
\begin{lemma}\label{lem:weights-parametrizations}
Letting $a=\Phi_\ii(b)$ and $y=\Upsilon_\ii(b)$, one has
\begin{align}
\wt(b)
&=-\sum_{j=1}^N a_j\alpha_{i_j}^\vee
 =-\sum_{j=1}^N y_j\beta_j^\vee
 &&\text{if }b\in B(\infty),\label{eq:weights-infinity}\\
\wt(b)
&=\lambda-\sum_{j=1}^N a_j\alpha_{i_j}^\vee
 =\lambda-\sum_{j=1}^N y_j\beta_j^\vee
 &&\text{if }b\in B(\lambda).
\label{eq:weights-finite}
\end{align}
\end{lemma}
\begin{proof}
Recall that the standard projection
$\pi_\lambda:B(\infty)\longrightarrow B(\lambda)\cup\{0\}$
restricts to a bijection from the surviving part of $B_w(\infty)$
onto $B_w(\lambda)$. We first show that

\vspace{0.3cm}
(a) $\wt(\pi_\lambda(b)) = \lambda+\wt(b)$
whenever $\pi_\lambda(b)\neq 0$.

\vspace{0.3cm}
Indeed, $\pi_\lambda$ sends $b_\infty$ to $b_\lambda$ and commutes
with the lowering operators $\widetilde f_i$; see
\cite[Theorem~2.3]{Fu22}. Thus, if
$b=\widetilde f_{i_r}\cdots\widetilde f_{i_1}b_\infty$ and
 $\pi_\lambda(b)\ne0$,
then $\pi_\lambda(b)
 =
 \widetilde f_{i_r}\cdots\widetilde f_{i_1}b_\lambda$.
Since $\wt(b_\infty)=0,
 \wt(b_\lambda)=\lambda,
 \wt(\widetilde f_i c)=\wt(c)-\alpha_i^\vee$
whenever $\widetilde f_i c\neq0$,
(a) follows.

By the definition of string coordinates, if $a=\Phi_\ii(b)$,
then
$b=\widetilde f_{i_1}^{\,a_1}\cdots \widetilde f_{i_N}^{\,a_N}b_\infty$ for $b\in B(\infty)$,
and $b=\widetilde f_{i_1}^{\,a_1}\cdots \widetilde f_{i_N}^{\,a_N}b_\lambda$
for $b\in B(\lambda)$.
This proves the string-coordinate equalities in
\eqref{eq:weights-infinity} and \eqref{eq:weights-finite}.

We next derive the finite-level formula for Lusztig coordinates
using $\Omega_{\ii,\lambda}$. Let $x=\Omega_{\ii,\lambda}(\Phi_\ii(b))
=\Upsilon_{\ii^*}(\eta_\lambda b)$,
and write $\beta_k^*=-w_0\beta_k,
 (\beta_k^*)^\vee=-w_0\beta_k^\vee$.

We claim that

\vspace{0.3cm}
(b) $\sum_{k=1}^N x_k\beta_k^\vee = \wt(b)-w_0\lambda, \sum_{k=1}^N x_k(\beta_k^*)^\vee = \lambda-w_0\wt(b)$.

\vspace{0.3cm}
Put $p_0=1$ and $p_j=s_{i_1}\cdots s_{i_j}$. For
$\gamma\in X_*(T)_{\mathbb R}$, we have
\[
 p_{j-1}\gamma-p_j\gamma
 =
 \langle\gamma,\alpha_{i_j}\rangle\beta_j^\vee.
\]
Summing over $1\leq j\leq N$ with $\gamma=\lambda$, and over
$j<k$ with $\gamma=\alpha_{i_k}^\vee$, respectively, yields
\[
 \sum_{j=1}^N
 \langle\lambda,\alpha_{i_j}\rangle\beta_j^\vee
 =
 \lambda-w_0\lambda, 
\text{ and }
 \beta_k^\vee+
 \sum_{j<k}
 \langle\alpha_{i_k}^\vee,\alpha_{i_j}\rangle\beta_j^\vee
 =
 \alpha_{i_k}^\vee.
\]

By the definition of $\Omega_{\ii,\lambda}$, we have
\[
 x_j
 =
 \langle\lambda,\alpha_{i_j}\rangle-a_j
 -
 \sum_{k>j}
 \langle\alpha_{i_k}^\vee,\alpha_{i_j}\rangle a_k.
\]
Therefore, collecting the coefficients of each $a_k$, we obtain
\[
 \begin{aligned}
 \sum_{j=1}^{N}x_j\beta_j^\vee
 &=
 \sum_{j=1}^{N}
 \left(
 \langle\lambda,\alpha_{i_j}\rangle-a_j
 -\sum_{k>j}
  \langle\alpha_{i_k}^\vee,\alpha_{i_j}\rangle a_k
 \right)\beta_j^\vee\\
 &=
 \lambda-w_0\lambda
 -\sum_{k=1}^{N}a_k
 \left(
 \beta_k^\vee+
 \sum_{j<k}
 \langle\alpha_{i_k}^\vee,\alpha_{i_j}\rangle\beta_j^\vee
 \right)\\
 &=
 \lambda-w_0\lambda
 -\sum_{k=1}^N a_k\alpha_{i_k}^\vee\\
 &=
 \wt(b)-w_0\lambda.
 \end{aligned}
\]
Applying $-w_0$ proves the second identity in~(b).

Since $\wt(\eta_\lambda b)=w_0\wt(b)$, the second identity
in~(b) gives $\wt(\eta_\lambda b) =
 \lambda-\sum_{j=1}^N x_j(\beta_j^*)^\vee$.
Apply this identity with $(\ii,b)$ replaced by
$(\ii^*,\eta_\lambda b)$. Since $\Omega_{\ii^*,\lambda}
 \bigl(\Phi_{\ii^*}(\eta_\lambda b)\bigr)=
 \Upsilon_{(\ii^*)^*}
 \bigl(\eta_\lambda(\eta_\lambda b)\bigr)=\Upsilon_\ii(b)=y$,
we obtain $\wt(b) =
 \lambda-\sum_{j=1}^N y_j\beta_j^\vee$.
This proves the Lusztig-coordinate equality in
\eqref{eq:weights-finite}.

Finally, let $b\in B(\infty)$ and choose $\lambda$ sufficiently
dominant so that $\pi_\lambda(b)\ne0$; this is possible, e.g., via Lemma~\ref{lem:demazure-stability}. By the definition of the
Lusztig parametrization on $B(\lambda)$ using the unique surviving
lift, we have $\Upsilon_\ii(\pi_\lambda(b))=\Upsilon_\ii(b)$.

The finite-level formula and (a) now give $\lambda+\wt(b) = \wt(\pi_\lambda(b)) =
 \lambda-\sum_{j=1}^N y_j\beta_j^\vee$, giving the Lusztig-coordinate equality in
\eqref{eq:weights-infinity}.
\end{proof}

We now conclude with the connection of longest DBG paths and stable Demazure multiplicity.

\begin{proof}[Proof of Theorem~\ref{thm:A}]
Given $b\in B_x(\infty)$, set $y=\Upsilon_{\ii}(b)$.
Theorem ~\ref{thm:fujita} places $y$ in $\bigcup_{\kk\in R(\ii,xw_0)}E_\kk \cap\Z^N,$ while Lemma~\ref{lem:complement} and
Lemma~\ref{lem:partition} give a unique $p\in\bar{\CP}$ with $y\in H_p$.
Assign the label $y_j$ to the edge at $j\in I(p)$. Integer points of
$H_p$ have nonnegative labels on upward edges and positive labels on
downward edges, exactly the required admissibility conditions.

Conversely, given such a labeled path, put its labels in the coordinates
indexed by $I(p)$ and put zero in the complementary coordinates.
The resulting vector lies in $H_p\cap\Z^N$. By Theorem~\ref{thm:fujita}
and the bijectivity of $\Upsilon_{\ii}$, it is the parameter vector of
a unique element of $B_x(\infty)$. The disjoint decomposition ensures
that these constructions are inverse.

Finally, the weight identity is
$-\wt(b)=\sum_{j=1}^N y_j\beta_j^\vee
=\sum_{j\in I(p)}y_j\beta_j^\vee=\wt(p,m)$.
Every fixed-weight fibre is finite: there are finitely many paths, and
the positive height of each $\beta_j^\vee$ bounds each label by the height of
$\delta$. Taking cardinalities gives the assertion on $m^{\text{st}}_x(\delta)$.
\end{proof}
For $p(q)\in \BZ[q], n\in \BZ_{\geq 0}$, we denote by $[p(q)]_{q^n}$ the coefficient of $q^n$ in $p(q)$. We can now give a character formula for $B_x(\infty)$.
\begin{corollary}\label{prop:path-character}
Fix a reflection order and let
\[
P_{x,\delta}(q)
:=
\sum_{\substack{p:1\rightarrow x\ \text{increasing and labeled path}\\
                 \operatorname{wt}(p)=\delta}}
q^{\ell(p)}.
\]

Then, as an identity of formal characters,
\[
\operatorname{ch}B_x(\infty)=\sum_{\delta\in Q^\vee_+}
[P_{x,\delta}(q)]_{q^{\ell(x)}}e^{-\delta},
\]
\end{corollary}
\begin{proof}
Every increasing path from $1$ to $x$ has length at most $\ell(x)$,
and the paths of length $\ell(x)$ are precisely the longest increasing
paths. Thus, for every $\delta\in Q^\vee_+$,
\[
[P_{x,\delta}(q)]_{q^{\ell(x)}}
=
\#\bigl\{
(p,m)\in\CP^\prec_{\max}(1\rightarrow x):
\operatorname{wt}(p,m)=\delta
\bigr\}.
\]
By Theorem~\ref{thm:A}, the latter set is in bijection with
$B_x(\infty)_{-\delta}$: explicitly, the bijection sends
$b\in B_x(\infty)$ to a path of weight $-\operatorname{wt}(b)$.
Consequently, $[P_{x,\delta}(q)]_{q^{\ell(x)}}
=
\#B_x(\infty)_{-\delta}$.
Multiplying by $e^{-\delta}$ and adding over $\delta\in Q^\vee_+$
gives
\[
\sum_{\delta\in Q^\vee_+}
[P_{x,\delta}(q)]_{q^{\ell(x)}}e^{-\delta}
=
\sum_{\delta\in Q^\vee_+}
\#B_x(\infty)_{-\delta}\,e^{-\delta}
=
\operatorname{ch}B_x(\infty).
\]
\end{proof}

\begin{lemma}\label{lem:demazure-stability}
Let $x\in W$ and $\delta=\sum_r d_r\alpha_r^\vee\in Q^\vee_+$.
Suppose that $\lambda$ satisfies $\langle\lambda,\alpha_r\rangle\ge 2d_r
~\text{for }\alpha_r\in \Delta$.
Then
\[
\dim V_x(\lambda)_{\lambda-\delta}
=
\#\{b\in B_x(\infty):\operatorname{wt}(b)=-\delta\}.
\]
In particular, this weight multiplicity is independent of $\lambda$
throughout the indicated range.
\end{lemma}

\begin{proof}
For $b\in B_x(\infty)$ with $\operatorname{wt}(b)=-\delta$,
write $a=\Phi_{\mathbf i}(b)$. Then $a_k\ge0$ and
\[
\delta=\sum_{k=1}^N a_k\alpha_{i_k}^\vee,
\qquad
d_r=\sum_{k:\,i_k=r}a_k.
\]

Since off-diagonal Cartan integers are nonpositive, for every $j$,
\[
\begin{aligned}
a_j+\sum_{k>j}a_k
\langle\alpha_{i_k}^\vee,\alpha_{i_j}\rangle
&\le a_j+2\sum_{\substack{k>j\\i_k=i_j}}a_k\le 2d_{i_j}\le \langle\lambda,\alpha_{i_j}\rangle.
\end{aligned}
\]
Thus $a$ satisfies all the highest-weight inequalities defining
$\Delta_{\mathbf i}(\lambda)$ inside the string cone, see Eq.~\ref{eq:string}.
Hence $\pi_\lambda(b)\ne0$. Consequently, every element of $B_x(\infty)_{-\delta}$ survives,
and the Demazure-compatible projection gives a bijection $B_x(\infty)_{-\delta}
\xrightarrow{\;\sim\;}
B_x(\lambda)_{\lambda-\delta}$.
Taking cardinalities yields
\[
\dim V_x(\lambda)_{\lambda-\delta}
=m_x^{\mathrm{st}}(\delta).
\]
\end{proof}

\begin{proof}[Proof of Theorem~\ref{thm:B}]
    Set $\delta=\lambda-\nu(b)=\sum d_i\alpha_i^\vee$ and note that under our depth hypothesis $\<\lambda,\alpha_i\> \geq 3\sum d_i \geq 2d_i$ for all $i$. The result now follows from Lemma~\ref{lem:demazure-stability}, Theorem~\ref{thm:DBGcount} and Theorem~\ref{thm:A}.
\end{proof}

\section{Relation between irreducible components of $X_w(b)$ and $X_{at^\gamma}(b)$}

\subsection{Admissible paths and associated strata}

In this section, we establish a decomposition of $X_w(b)$ into certain locally closed subschemes by refining techniques from 
\cite{GH10}. This is in the spirit of a similar decomposition using the Deligne-Lusztig reduction, but with certain specific advantages as follows. The leaves for Deligne-Lusztig reduction tree for $w$ are certain elements of $\widetilde{W}$ that are minimal length in their conjugacy class, and their cardinality grows as $w$ gets bigger; on the other hand, our decomposition scheme ends at certain elements of the form encountered in the previous section, and the cardinality of such elements is bounded independently of $\lambda$, the translation part of $w$, in a precise sense as we see below.

Recall that we have the map of relative position $\inv: \text{Fl}\times \text{Fl} \rightarrow \widetilde W$ defined by $g_1^{-1}g_2 \in I~\inv(g_1I, g_2I)~I$. Following \cite[\S 2.4]{GH10}, we adopt a convenient notational device to describe varieties of tuples of elements in $
\mathrm{Fl}$: each braced arrow diagram denotes the locally closed
variety of tuples (inside the appropriate product of affine flag varieties) satisfying the indicated relative positions. We freely move between specifying such varieties either by arrow diagrams or by listing out the relevant relative positions (especially if drawing the diagram gets visually complicated). We now collect some standard facts below. First, we have an elementary left-handed version of the rank-one phenomenon from \cite[\S 2.4]{GH10}.
\begin{lemma}\label{lem:left-handed}
Let $s\in\widetilde{\mathbb S}$ and $w\in\widetilde W$.
\begin{enumerate}[(i)]
\item If $sw>w$, then $\pi: (g,g'',g')\longmapsto(g,g')$ 
gives an isomorphism $\left\{g\xrightarrow{s}g''\xrightarrow{w}g'\right\}
\longrightarrow
\left\{g\xrightarrow{sw}g'\right\}$.
\item If $sw<w$, there is a disjoint decomposition
\[
\left\{g\xrightarrow{s}g''\xrightarrow{w}g'\right\} =
\left\{
g\xrightarrow{s}g''\xrightarrow{w}g'
\;\middle|\;
\operatorname{inv}(g,g')=sw
\right\} \sqcup
\left\{
g\xrightarrow{s}g''\xrightarrow{w}g'
\;\middle|\;
\operatorname{inv}(g,g')=w
\right\}.
\]
Call the summands $X_1$, resp. $X_2$. Then $X_1$ is closed, $X_2$ is open. Further, 
the projections $\pi: X_1 \rightarrow \left\{g\xrightarrow{sw}g'\right\}$ and $\pi: X_2 \rightarrow \left\{g\xrightarrow{w}g'\right\}$ are surjective
Zariski-locally trivial bundles with fibres
$\mathbb A^1$ and $\mathbb G_m$, respectively.
\end{enumerate}
 
\end{lemma}

Fix $w_1,w_2\in\widetilde W$ and set $\CC_{w_2,w_1}(b):=\{g_1\xrightarrow {w_2} g_2\xrightarrow {w_1} b\sigma(g_1)\}.$ Note that if $\ell(w_1w_2)=\ell(w_1)+\ell(w_2)$, we have $X_{w_1w_2}(b) \simeq \{g\xrightarrow {w_1} g_1\xrightarrow {w_2} b\sigma(g)\}$. The following is then a straightforward observation, see \cite[\S 5.4]{He21b}.

\begin{lemma}\label{lem:rotation}
If $w=w_1w_2$ and $\ell(w)=\ell(w_1)+\ell(w_2)$, then the map $(g,g_1) \mapsto (g_1,b\s(g))$ is a $J_b(F)$-equivariant universal homeomorphism $h:X_w(b)\longrightarrow  \CC_{w_2,w_1}(b)$

\end{lemma}
Let $m=\ell(w_2)$ and fix throughout this section a reduced expression $w_2=r_1\cdots r_m\tau$, with $r_i\in\widetilde{\mathbb S}$ and $\tau\in\Omega$. A length-zero element gives a unique algebraic step: its representative normalizes $I$. Thus $A_m\xrightarrow\tau H\xrightarrow{w_1} B$ is equivalent to $A_m\xrightarrow{\tau w_1}B$. Combining with Lemma~\ref{lem:left-handed}[(i)], the $w_2$-step then has a unique subdivision resulting in
\begin{equation}\label{eq:subdivision}
 \CC_{w_2,w_1}(b) \simeq \{g_1=a_0\xrightarrow{r_1}a_1\xrightarrow{r_2}\cdots
 \xrightarrow{r_m}a_m\xrightarrow\tau g_2 \xrightarrow{w_1} b\sigma(g_1)\}.
\end{equation}
We may thus obtain a stratification of $\CC_{w_2,w_1}(b)$ by \emph{adding more arrows}, i.e., by prescribing relative positions between the intermediate flags $a_i$ and the final flag $b\sigma(g)$. In view of Lemma~\ref{lem:left-handed}, this is formalized in the following definition.

\begin{definition}\label{def:path}
An \emph{admissible path} is a sequence $\CP=(v_0,\ldots,v_m)\in \widetilde W^{m+1}$, such that $v_m=\tau w_1$ and for $i=1,\cdots,m$,
\begin{equation}\label{eq:path-rule}
 v_{i-1}=\begin{cases}
 r_i v_i,&r_i v_i>v_i,\\
 r_i v_i\ \text{or}\ v_i,&r_i v_i<v_i.
 \end{cases}
\end{equation}
We define its \emph{endpoint} to be $v_{\CP}=v_0$. Let $n_A(\CP)$ count the indices with $r_i v_i<v_i$ and $v_{i-1}=r_i v_i$, and let $n_G(\CP)$ count those with $v_{i-1}=v_i$. Put $\delta_{\CP}=n_A(\CP)+n_G(\CP)$.
\end{definition}
We now give these paths a geometric interpretation.
\begin{proposition}\label{prop:strata}
For each admissible path $\CP$, define
\[
\CS_{\CP}(b)
=
\left\{
(g_1,g_2)\in\CC_{w_2,w_1}(b)
\;\middle|\;
\inv(a_i,b\sigma(g_1))=v_i
\text{ for }0\le i\le m
\right\},
\]
where $a_0=g_1,a_1,\ldots,a_m$ are the intermediate flags
determined by the unique subdivision in Eq.~\ref{eq:subdivision}. Then these loci give a decomposition of $\CC_{w_2,w_1}(b)$ into finitely many $J_b(F)$-stable locally closed subvarieties.

The restriction of the projection $(g_1,g_2) \mapsto g_1$ defines
a morphism $\pi_{\CP}: S_{\CP}(b)\rightarrow X_{v_{\CP}}(b)$, that is a surjective tower
of Zariski-locally trivial $\BA^1$- and $\BG_m$-bundles (including possible isomorphism steps), of relative dimension $\delta_{\CP}$. Moreover,
\begin{equation}\label{eq:length-count}
 \ell(w)-\ell(v_{\CP})=2n_A(\CP)+n_G(\CP).
\end{equation}
\end{proposition}

\begin{proof}
Every point in $\CC_{w_2,w_1}(b)$ determines a sequence of relative positions $v_i$ via Eq.~\ref{eq:subdivision}. Lemma~\ref{lem:left-handed} gives precisely Eq.~\ref{eq:path-rule}, so these loci exhaust $\CC_{w_2,w_1}(b)$ and are disjoint. 

For a fixed path $\CP$, let $Z_i$ be the reduced locally closed subscheme
\[
Z_i=
\left\{
(a_0,\ldots,a_i)\in\mathrm{Fl}^{i+1}
\;\middle|\;
\begin{array}{ll}
\inv(a_{j-1},a_j)=r_j & \text{for }1\le j\le i,\\
\inv(a_j,b\sigma(a_0))=v_j & \text{for }0\le j\le i
\end{array}
\right\}.
\]

Then, $Z_0=X_{v_0}(b)$ and $Z_m\simeq \CS_{\CP}(b)$. For $1\le i\le m$, let $f_i\colon Z_i\rightarrow Z_{i-1},
(a_0,\ldots,a_i)\mapsto(a_0,\ldots,a_{i-1})$.
The fibre over a fixed $(a_0,\ldots,a_{i-1})$ then is
\[
 \{h:\inv(a_{i-1},h)=r_i,\ \inv(h,b\sigma(a_0))=v_i\}.
\]

By Lemma~\ref{lem:left-handed} and admissibility of $\CP$, the map $f_i$ is
\[
\begin{cases}
\text{an isomorphism},
& r_iv_i>v_i,\\
\text{a surjective Zariski-locally trivial $\BA^1$-bundle},
& r_iv_i<v_i,\quad v_{i-1}=r_iv_i,\\
\text{a surjective Zariski-locally trivial $\BG_m$-bundle},
& r_iv_i<v_i,\quad v_{i-1}=v_i.
\end{cases}
\]
Consequently,
\[
\pi_{\CP}\colon
\CS_{\CP}(b)\simeq Z_m
\xrightarrow{f_m}Z_{m-1}\longrightarrow\cdots
\xrightarrow{f_1}Z_0=X_{v_0}(b)
\]
is the asserted surjective tower, of relative dimension
$n_A(\CP)+n_G(\CP)=\delta_{\CP}$.

Since every index belongs to exactly one of these three types, $m=n_++n_A+n_G$.
Summing the length changes gives
\[
\ell(v_0)-\ell(v_m)
=\sum_{i=1}^m\bigl(\ell(v_{i-1})-\ell(v_i)\bigr)
=n_+\cdot 1 + n_A\cdot(-1) + n_G \cdot 0= n_+-n_A.
\]
Therefore
\[
\begin{aligned}
\ell(w)-\ell(v_0)
&=m-\bigl(\ell(v_0)-\ell(w_1)\bigr)\\
&=(n_++n_A+n_G)-(n_+-n_A)\\
&=2n_A+n_G,
\end{aligned}
\]
which proves \eqref{eq:length-count}.

Finally, $b\sigma(jg)=j b\sigma(g)$ for $j\in J_b(F)$, so the conditions and every forgetting map are equivariant.
\end{proof}

\subsection{Special endpoint elements}\label{subsec:cordial ends}
For $w\in\widetilde{W}$ write it in
standard form $w=xt^{\lambda}y; x,y\in W, t^{\lambda}y\in{}^{\BS}\tW$. Let $J=\{\alpha_i:y^{-1}\alpha_i <0\} \subset \Delta$, and choose a dominant cocharacter $\mu$ such that $\langle\mu,\alpha_i\rangle\ge1$ for $\alpha_i\in J$. This condition ensures that $t^\mu y\in{}^{\mathbb S}\widetilde W$ and
thus $\ell(t^\mu y)=\ell(t^\mu)-\ell(y)$. Fix this choice independently of $\lambda$.

Regarding the choice of $\mu$, we may take $\mu=\rho_J^\vee:=\sum_{\alpha_i\in J}\omega_i^\vee$ when it belongs to $X_*(T)$. In every root datum, $\mu=2\rho^\vee:=\sum\limits_{\a \in \Phi^+} \alpha^\vee$ is an available choice. Set $M_\mu=\max\limits_{\alpha\in\Phi^+}\langle\mu,\alpha\rangle=\langle \mu, \theta \rangle$. Write $\theta=\sum n_i\alpha_i$, then $\sum n_i=h-1$.  If $\mu=\rho^\vee_J$, $M_\mu = \sum\limits_{i\in J} n_i \leq h-1$, where $h$ is the Coxeter number for $\Phi$; if $\mu=2\rho^\vee$, then $M_\mu\leq 2h-2$ by the same token.

Assume further that $\kappa:=\lambda-\mu$ is dominant; the regularity hypothesis imposed in Lemma~\ref{prop:regularity}
automatically ensures this. Set $w_1=xt^\kappa, w_2=t^\mu y$. Then $w=w_1w_2$ and $\ell(w)=\ell(w_1)+\ell(w_2)$. For an admissible path $\CP$, define
\begin{equation}\label{eq:path-epsilon}
 \varepsilon_i=\begin{cases}
 1,&v_{i-1}=r_iv_i,\\
 0 ,&v_{i-1}=v_i.
 \end{cases}
\end{equation}

Then $v_{\CP}=u_{\CP}xt^\kappa$, where
 $u_{\CP}:=r_1^{\varepsilon_1}\cdots r_m^{\varepsilon_m}\tau\le w_2$. Writing $u_{\CP}=c_{\CP}t^{\beta_{\CP}}$, with $c_{\CP}\in W$, we then have
 \begin{equation}\label{eq:a-gamma}
 v_{\CP}=a_{\CP}t^{\gamma_{\CP}}, \text{ where }
 a_{\CP}=c_{\CP}x,
 \gamma_{\CP}=\lambda-\mu+x^{-1}\beta_{\CP}.
\end{equation}

\begin{lemma}\label{prop:regularity}
If $\langle\lambda,\alpha_i\rangle\ge\langle\mu,\alpha_i\rangle+M_\mu$ for every simple root $\alpha_i$, then every $\gamma_{\CP}$ is dominant; consequently, every endpoint $a_{\CP}t^{\gamma_{\CP}}$ is cordial.
\end{lemma}
\begin{proof}
Since $t^\mu=w_2y^{-1}$ and
$\ell(t^\mu)=\ell(w_2)+\ell(y^{-1})$, we have $w_2\leq t^\mu$.
Consequently, $u_{\CP}\leq t^\mu$. Compatibility of the Bruhat
order on $\widetilde W$ with the dominance order under the
projection $Wt^\zeta W\mapsto\zeta^+$ then gives $\beta_{\CP}^{+}\leq\mu^{+}=\mu$,
where $\zeta^+$ denotes the dominant representative in the orbit $W\zeta$.
By the interpretation of dominance order in terms of convex hulls
\cite[Lemma~12.14]{AB83}, it follows that
$x^{-1}\beta_{\CP}\in\operatorname{Conv}(W\mu)$,
the convex hull of the $W$-orbit of $\mu$. We now claim:

\vspace{0.3cm}
(a) every $\zeta\in\text{Conv}(W\mu)$ satisfies $\langle\zeta,\alpha_i\rangle\geq -M_\mu$.

\vspace{0.3cm}
It suffices to check the vertices $z\mu$, where the pairing is $\langle\mu,z^{-1}\alpha_i\rangle\ge \langle\mu,-\theta\rangle = -M_\mu$. Therefore
\[
 \langle\gamma_{\CP},\alpha_i\rangle
 \geq\langle\lambda,\alpha_i\rangle-\langle\mu,\alpha_i\rangle-M_\mu\geq 0.
\]
\end{proof}
\begin{remark}\label{rmk:depth}
 The last estimate in the proof above shows that
\[
\operatorname{depth}(\gamma_{\CP})
\geq
\operatorname{depth}(\lambda)
-\max_i\langle\mu,\alpha_i\rangle-M_\mu.
\]

If we may choose a dominant cocharacter $\mu=\rho_J^\vee$, then $M_\mu\leq h-1$ and hence
$\operatorname{depth}(\gamma_{\CP})
\geq\operatorname{depth}(\lambda)-h$.
In particular, $\operatorname{depth}(\lambda)\geq h+A$ implies
$\operatorname{depth}(\gamma_{\CP})\geq A$. For an arbitrary root datum, we may always take $\mu=2\rho^\vee$.
Then $M_\mu\leq2(h-1)$, and the same argument gives instead $\operatorname{depth}(\gamma_{\CP})
\geq\operatorname{depth}(\lambda)-2h$.
In particular, the hypothesis of Lemma~\ref{prop:regularity} is satisfied
if $\operatorname{depth}(\lambda)\geq h$ for the first choice of
$\mu$, and if $\operatorname{depth}(\lambda)\geq2h$ for the choice
$\mu=2\rho^\vee$.
\end{remark}

\begin{theorem}\label{thm:main}
Assume the regularity condition of Lemma~\ref{prop:regularity}. Let $\mathscr{P}=\mathscr{P}(w)$ be the collection of admissible paths starting at $\tau x t^\kappa$. Then we have the following assertions.
\begin{enumerate}[(i)]
\item The nonemptiness set is
\begin{equation}\label{eq:B-union}
 B(G)_w=\bigcup_{\CP\in\mathscr{P}}B(G)_{a_{\CP}t^{\gamma_{\CP}}}.
\end{equation}
\item If $X_w(b)\neq\emptyset$, then
\begin{equation}\label{eq:dimension-max}
 \dim X_w(b)=\max_{\substack{\CP\in\mathscr{P}:\\ [b]\in B(G)_{v_{\CP}}}}
 \bigl(d_{v_{\CP}}(b)+\delta_{\CP}\bigr).
\end{equation}
\item Put $d=\dim X_w(b)$ and define
\[
 \mathscr{P}(w)_b^{\mathrm{top}}=
 \{\CP\in\mathscr{P}:X_{v_{\CP}}(b)\neq\emptyset,
 \ d_{v_{\CP}}(b)+\delta_{\CP}=d\}.
\]
There is a $J_b(F)$-equivariant bijection
\begin{equation}\label{eq:main-components}
 \bigsqcup_{\CP\in\mathscr{P}(w)_b^{\mathrm{top}}}
 \Sigma\bigl(X_{a_{\CP}t^{\gamma_{\CP}}}(b)\bigr)
 \xrightarrow{\sim}\Sigma^{\text{top}}\bigl(X_w(b)\bigr).
\end{equation}
\end{enumerate}

\end{theorem}
Note that the number of paths is at most $2^{\ell(t^\mu y)}$, independently of $\lambda$. In \eqref{eq:B-union} repeated endpoints can be removed; in \eqref{eq:main-components} paths are retained with their multiplicities.

\begin{proof}
Lemma~\ref{lem:rotation} and Proposition~\ref{prop:strata} give a
finite partition
\[
\CC_{w_2,w_1}(b)=\bigsqcup_{\CP\in\mathscr P} S_{\CP}(b)
\]
and, for each $\CP$, a surjective tower
$\pi_{\CP}:S_{\CP}(b)\longrightarrow X_{v_{\CP}}(b)$
of relative dimension $\delta_{\CP}$. Thus $S_{\CP}(b)$ is
nonempty if and only if $X_{v_{\CP}}(b)$ is nonempty. Since
$X_w(b)$ and $\CC_{w_2,w_1}(b)$ are universally homeomorphic,
this proves \eqref{eq:B-union}.

Suppose now that $X_w(b)$ is nonempty. By
Lemma~\ref{prop:regularity}, every $v_{\CP}$ is cordial.
Hence, whenever $X_{v_{\CP}}(b)$ is nonempty, it is equidimensional
of dimension $d_{v_{\CP}}(b)$. The bundle tower then shows that
$S_{\CP}(b)$ is equidimensional of dimension
$d_{v_{\CP}}(b)+\delta_{\CP}$. Since the partition is finite,
\[
\dim X_w(b)
=
\max_{\substack{\CP\in\mathscr P:\\
                 X_{v_{\CP}}(b)\ne\varnothing}}
\bigl(d_{v_{\CP}}(b)+\delta_{\CP}\bigr),
\]
proving \eqref{eq:dimension-max}.

It remains to compare irreducible components. For each $\CP$, the
fibres of $\pi_{\CP}$ are obtained by successive
$\mathbb A^1$-, $\mathbb G_m$-, and isomorphism steps. Since the
first two have irreducible fibres and are Zariski-locally trivial,
inverse image gives a bijection
$\Sigma\bigl(X_{v_{\CP}}(b)\bigr)
\xrightarrow{\sim}
\Sigma\bigl(S_{\CP}(b)\bigr)$.

Put $d=\dim X_w(b)$. Each $d$-dimensional component of
$\CC_{w_2,w_1}(b)$ contains a dense open subset of a unique
stratum $S_{\CP}(b)$ of dimension $d$. Conversely, the closure of
an irreducible component of such a stratum is a $d$-dimensional
component of $\CC_{w_2,w_1}(b)$. Consequently, taking closures gives
a bijection
\[
\bigsqcup_{\CP\in\mathscr P(w)_b^{\mathrm{top}}}
\Sigma\bigl(S_{\CP}(b)\bigr)
\xrightarrow{\sim}
\Sigma^{\mathrm{top}}\bigl(\CC_{w_2,w_1}(b)\bigr).
\]
Composing with the preceding bijections and transporting components
through the universal homeomorphism
$h:X_w(b)\to\CC_{w_2,w_1}(b)$ proves
\eqref{eq:main-components}. Explicitly, a component
$C\subseteq X_{v_{\CP}}(b)$ is sent to
$h^{-1}\!\left(
\overline{\pi_{\CP}^{-1}(C)}^{\,\CC_{w_2,w_1}(b)}
\right)$.
All the strata, maps $\pi_{\CP}$, and $h$ are $J_b(F)$-equivariant,
so the resulting bijection is equivariant.
\end{proof}

\begin{corollary}\label{cor:orbits}
The bijection in \eqref{eq:main-components} gives
\[
 \#~ J_b(F)\backslash\Sigma^{top}(X_w(b))=\sum_{\CP\in\mathscr{P}(w)_b^{\mathrm{top}}}\#~J_b(F)\backslash\Sigma(X_{a_{\CP}t^{\gamma_{\CP}}}(b)).
\]

\end{corollary}
\begin{proof}
Each stratum is $J_b(F)$-stable. Taking orbits of the equivariant disjoint union therefore gives the stated disjoint union of orbit sets. 
\end{proof}

\begin{proof}[Proof of Theorem~\ref{thm:C-strong}] In view of Lemma~\ref{prop:regularity} and Remark~\ref{rmk:depth}, the assertion in (i) follows from Theorem~\ref{thm:main}[(iii)].

For $\CP\in\mathscr P(w)_b^{\mathrm{top}}$, we have $X_{v_{\CP}}(b)\neq \emptyset$, thus $\nu(b)\leq\gamma_{\CP}$. Moreover, the convex-hull argument in
Lemma~\ref{prop:regularity} gives $\gamma_{\CP}\leq\lambda$ as follows. We have
$x^{-1}\beta_{\CP}\in\operatorname{Conv}(W\mu)$, and since $\mu$ is
dominant, every point of this convex hull is smaller than $\mu$ in the dominance order, giving $\lambda-\gamma_{\CP}
=\mu-x^{-1}\beta_{\CP}
\in\sum_i\mathbb R_{\geq0}\alpha_i^\vee$,
so $\gamma_{\CP}\leq\lambda$.

Thus
$\nu(b)\leq\gamma_{\CP}\leq\lambda$, ensuring that 
$\langle\gamma_{\CP}-\nu(b),\rho\rangle
\leq\langle\lambda-\nu(b),\rho\rangle$. The depth hypothesis on $\lambda$ and Remark~\ref{rmk:depth} therefore imply
\[
\operatorname{depth}(\gamma_{\CP})
\geq 3\langle\lambda-\nu(b),\rho\rangle
\geq 3\langle\gamma_{\CP}-\nu(b),\rho\rangle.
\]
We can consequently apply Theorem~\ref{thm:B} to each ADLV at the nonempty endpoints;
the formula follows from Corollary~\ref{cor:orbits}.
\end{proof}

\begin{example}
We give an example in $G=\mathrm{PGL}_3$. Write $\omega_1^\vee,\omega_2^\vee$ for the fundamental coweights, and put $\rho^\vee=\omega_1^\vee+\omega_2^\vee,
w_0=s_1s_2s_1, s_0=t^{\rho^\vee}w_0$. Take $y=s_1s_2,
\mu=2\omega_1^\vee, \kappa=3\rho^\vee,
\lambda=\mu+\kappa=5\omega_1^\vee+3\omega_2^\vee$.
We choose $\mu$ strictly larger than $\rho_J^\vee$ in order to obtain a
two-step decomposition. Set $w=t^\lambda y$ and
\[
w_2=t^\mu y,\qquad w_1=t^\kappa,\qquad
\tau=t^{\omega_2^\vee}s_2s_1.
\]
Then $\ell(\tau)=0$ and $w_2=s_0s_2\tau, r:=\tau t^\kappa
=s_2s_1t^{2\omega_1^\vee+3\omega_2^\vee}$.
It can be directly seen that the regularity condition of Lemma~\ref{prop:regularity} holds.

A point of $\CC_{w_2,w_1}(b)$ has a unique subdivision
\[
g_1=a_0\xrightarrow{s_0}a_1\xrightarrow{s_2}a_2
\xrightarrow{\tau}g_2\xrightarrow{t^\kappa}b\sigma(g_1),
\]
and $\inv(a_2,b\sigma(g_1))=r$. The relevant lengths are $\ell(r)=12, \ell(s_2r)=11, \ell(s_0s_2r)=10, \ell(s_0r)=13$. Thus, there are three admissible paths. Their strata and endpoints are as follows; the bundle steps are listed in the order of
forgetting $a_2$, then $a_1$.
\[
\begin{array}{c|c|c|c|c}
\text{Stratum}
&\inv(a_1,b\sigma(g_1))
&\inv(g_1,b\sigma(g_1))
&\text{Bundle steps}
&\delta_{\CP}\\ \hline
\CS_1(b)&s_2r&s_0s_2r&\mathbb A^1,\mathbb A^1&2\\
\CS_2(b)&s_2r&s_2r&\mathbb A^1,\mathbb G_m&2\\
\CS_3(b)&r&s_0r&\mathbb G_m,\text{isomorphism}&1
\end{array}
\]

For example, the projection of $S_2(b)$ to $X_{s_2r}(b)$ factors as
\[
S_2(b)
\longrightarrow
\left\{
g_1\xrightarrow{s_0}a_1\xrightarrow{s_2r}b\sigma(g_1):
g_1\in X_{s_2r}(b)
\right\}
\longrightarrow X_{s_2r}(b),
\]
where the first map is an $\mathbb A^1$-bundle and the second is a
$\mathbb G_m$-bundle.

Writing the endpoints as $v_i=a_i't^{\gamma_i}$, we obtain
$v_1=s_0s_2r =s_1s_2t^{3\omega_1^\vee+\omega_2^\vee};
v_2=s_2r=s_1t^{2\omega_1^\vee+3\omega_2^\vee};
v_3=s_0r=s_2t^{4\omega_1^\vee+2\omega_2^\vee}.$
All three translation parts are strictly dominant. In particular,
Theorem~\ref{thm:main}(i) gives
\[
B(G)_w=B(G)_{v_1}\cup B(G)_{v_2}\cup B(G)_{v_3}.
\]

Furthermore, $\ell(w)=14$, $\ell(\eta(w))=2$, and $\langle\gamma_i,\rho\rangle+\ell(a_i')+\delta_{\CP_i}=8
~\forall~ 1\leq i\leq3.$
The dimension formula for the cordial endpoints therefore shows that every nonempty stratum has dimension $d_w(b)$.
Consequently, whenever $X_w(b)$ is nonempty, it is equidimensional,
and Theorem~\ref{thm:main}(iii) gives a $J_b(F)$-equivariant
bijection
\[
\bigsqcup_{i=1}^{3}\Sigma\bigl(X_{v_i}(b)\bigr)
\xrightarrow{\sim}
\Sigma^{\mathrm{top}}\bigl(X_w(b)\bigr).
\]
\end{example}
\section{Some instances of multiplicity one phenomena}

In this section, we collect some instances where the left-hand side of the bijection in Eq.~\ref{eq:main-components} has a single term.
We start by identifying a particular summand in 
Eq.~\ref{eq:main-components} under the regularity condition of Lemma~\ref{prop:regularity}. Define
the path $\CP_0$ that chooses multiplication at every index: then $\CP_0$ is always admissible. Its endpoint is
\begin{equation}\label{eq:ordinary-endpoint}
 v_{\CP_0}=w_2w_1=yxt^{\gamma_0},\text{ where }
 \gamma_0=\lambda-\mu+(yx)^{-1}\mu.
\end{equation}
The path $\CP_0$ is in fact the only admissible path whose
endpoint has finite Weyl part $\eta(w)$. We will see this after Lemma~\ref{lem:stationary-finite}; using this endpoint, we obtain the following estimate.

\begin{corollary}\label{cor:ordinary-injection}
Assume the regularity condition of Lemma~\ref{prop:regularity} and $X_{yxt^{\gamma_0}}(b)\neq\emptyset$. Then $\dim X_w(b)=d_w(b)$ and the path $\CP_0$ gives a $J_b(F)$-equivariant injection
\[
 \Phi_0:\Sigma\bigl(X_{yxt^{\gamma_0}}(b)\bigr)
 \hookrightarrow\Sigma^{\text{top}}\bigl(X_w(b)\bigr).
\]
This injection is bijective if and only if for every other admissible path $\CP$ we have either $X_{v_{\CP}}(b)=\emptyset$ or $d_{v_{\CP}}(b)+\delta_{\CP}<d_w(b)$.
\end{corollary}
\begin{proof}
We note that $\CP_0$ has no \(\mathbb G_m\)-steps, its steps are \(\mathbb A^1\)-bundles or isomorphisms; since $n_G(\CP_0)=0$, we have $\delta_{\CP_0}=\frac{\ell(w)-\ell(v_{\CP_0})}{2}$ by Eq.~\ref{eq:length-count}. Since $\eta(w)=yx=\eta(v_{\CP_0})$, we get by direct calculation $\delta_{\CP_0}=d_w(b)-d_{v_{\CP_0}}(b)$. The \emph{ordinary stratum} $\CS_{\CP_0}(b)$ therefore has dimension $d_w(b)$, since the ADLV at the base satisfies $\dim X_{v_{\CP_0}}(b)=d_{v_{\CP_0}}(b)$ by cordiality of $v_{\CP_0}$. The general upper bound then forces $\dim X_w(b)=d_w(b)$. The component injection and its precise surjectivity criterion follow from \eqref{eq:main-components}.
\end{proof}

\begin{corollary}\label{cor:score}
Under the regularity condition of Lemma~\ref{prop:regularity}, we have for every admissible path $\CP$ 
\begin{equation*}
 d_w(b)-\{d_{v_{\CP}}(b)+\delta_{\CP}\}
 =\frac{\ell(yx)-\{n_G(\CP)+\ell(a_{\CP})\}}{2}\geq 0.
\end{equation*}
If $X_{v_{\CP_0}}(b)\neq \emptyset$, the stratum $\CS_{\CP}(b)$ contributes top components exactly when $X_{v_{\CP}}(b)$ is nonempty and $n_G(\CP)+\ell(a_{\CP})=\ell(yx)$.
\end{corollary}
\begin{proof}
The equality follows from the definition of virtual dimension and \eqref{eq:length-count}. Now, fix a path $\CP$ and choose $[b]=[\dot t^{\gamma_{\CP}}]$. Since $\gamma_{\CP}$ is dominant, this is the maximal $\sigma$-conjugacy class in
$B(G)_{v_{\CP}}=B(G)_{a_{\CP}t^{\gamma_{\CP}}}$; then its nonempty stratum $\CS_{\CP}(b)$ has dimension $d_{v_{\CP}}(b)+\delta_{\CP}$, bounded above by $\dim X_w(b)\leq d_w(b)$, proving the inequality. The final assertion uses Corollary~\ref{cor:ordinary-injection}.
\end{proof}
We now proceed to give a sufficient criterion under which the ordinary path $\CP_0$ is the only relevant one, i.e. $\CS_{\CP_0}(b)$ is the only top dimensional stratum in $X_w(b)$ as in Proposition~\ref{prop:strata}. Recall that we have fixed a reduced word $w_2=r_1\cdots r_m \tau$. Set $a_i:=r_1\cdots r_{i-1}$. The reduced expression of $w_2$ determines a sequence of crossings of hyperplanes to reach from $\mathfrak{a}$ to $w_2\mathfrak{a}$, and its $i$-th wall has the associated affine reflection $R_i:=a_ir_ia_i^{-1}$. 

For an admissible path $\CP=(v_0,\ldots,v_m)$ with
$v_m=\tau w_1$ as before, set $U_i:=(r_1\cdots r_i)v_i$ for $0\leq i\leq m,$
using the empty product when $i=0$; thus,
$U_m=w_2w_1$ and $U_0=v_{\CP}$.
In order to compare the paths $\CP_0$ and $\CP$, we need to understand the effect of \emph{stationary choices}, i.e. the steps in $\CP$ where $v_{i-1}=v_i$. This is recorded in the lemma below.
\begin{lemma}\label{lem:stationary-finite}
We have
\begin{equation}\label{eq:nested-labels}
 R_i=t^{k\beta^\vee}s_\beta,\text{ for certain }\beta\in\Phi^+,
 1\leq k\leq\langle\mu,\beta\rangle -\mathbf1_{\{y^{-1}\beta<0\}}.
\end{equation}
Writing $U_i=zt^{\xi}$, we further have
\[
U_{i-1}
=
\begin{cases}
zt^{\xi},&v_{i-1}=r_iv_i;\\
(s_\beta z)t^{\xi-kz^{-1}\beta^\vee},&v_{i-1}=v_i,
\end{cases}
\]
where in the second case $z^{-1}\beta<0$, equivalently $s_\beta z < z$.

\end{lemma}
\begin{proof}
Let $\mathfrak a$ be the base alcove. The element $w_2=t^\mu y$ is minimal in its left $W$-coset, hence the alcove $w_2\mathfrak{a}$ lies in the dominant chamber.
For $c\in\mathfrak a$ and $\beta\in\Phi^+$,
\[
\langle w_2c,\beta\rangle
=\langle\mu,\beta\rangle+\langle c,y^{-1}\beta\rangle.
\]
The second term lies in $(0,1)$ if $y^{-1}\beta>0$,
and in $(-1,0)$ otherwise.
Since $w_2\mathfrak a$ lies in the dominant chamber,
the hyperplanes separating $\mathfrak a$ and $w_2\mathfrak a$
are precisely
\[
H_{\beta,k}=\{v:\langle v,\beta\rangle=k\},
\qquad \beta\in\Phi^+,\quad 1\le k\le \langle\mu,\beta\rangle
-\mathbf1_{\{y^{-1}\beta<0\}}.
\]
The reflections $R_i$ enumerate these separating hyperplanes,
proving \eqref{eq:nested-labels}.

From the definition of $U_i$, we have
\[
U_{i-1}
=
\begin{cases}
U_i,&v_{i-1}=r_iv_i,\\
R_iU_i,&v_{i-1}=v_i.
\end{cases}
\]

It remains to establish the claim that $z^{-1}\beta<0$ in the second case. Let
$\alpha_{r_i}$ be the simple affine root corresponding to $r_i$.
By the admissibility rule, the second possibility above occurs precisely when
$r_i v_i<v_i$, equivalently when $v_i^{-1}\alpha_{r_i}<0$.
Since the positive affine root corresponding to $R_i$ is $(r_1\cdots r_{i-1})(\alpha_{r_i})
=k\delta_{\mathrm{aff}}-\beta$, in view of
the fact $U_i=(r_1\cdots r_{i-1})r_i v_i$ and
$r_i(\alpha_{r_i})=-\alpha_{r_i}$, the condition is then equivalent
to $U_i^{-1}(k\delta_{\mathrm{aff}}-\beta)>0.$ Using $U_i=zt^\xi$, we calculate
\[
U_i^{-1}(k\delta_{\mathrm{aff}}-\beta)
=
-z^{-1}\beta+
\bigl(k-\langle\xi,z^{-1}\beta\rangle\bigr)
\delta_{\mathrm{aff}}.
\tag{*}
\]

Both $U_i$ and $R_iU_i$ have the form $rw_1$, where $r$ is subword of $w_2=r_1\cdots r_m\tau$. In particular,
$r\leq w_2\leq t^\mu$. The convex-hull argument in the proof of Lemma~\ref{prop:regularity} therefore applies to both elements: their  translation parts $\xi$ and
$\xi'=\xi-kz^{-1}\beta^\vee$ are dominant. If $z^{-1}\beta<0$, then
dominance of $\xi$ makes the coefficient of $\delta_{\mathrm{aff}}$
in $(*)$ positive. If $z^{-1}\beta>0$, dominance of $\xi'$ gives
\[
\langle\xi,z^{-1}\beta\rangle
=
\langle\xi',z^{-1}\beta\rangle+2k
\geq 2k,
\]
so that coefficient is negative. Hence $(*)$ is positive exactly
when $z^{-1}\beta<0$.
\end{proof}
In particular, if $\CP\ne\CP_0$, then $\CP$ has at least one
stationary choice. Starting from $U_m=yxt^{\gamma_0}$ and accounting for such choices as in the proof above, Lemma~\ref{lem:stationary-finite} shows that $\ell(a_{\CP})<\ell(yx)$.
Thus, $\CP_0$ is the unique admissible path whose endpoint has
finite Weyl part $\eta(w)=yx$.
\begin{proposition}
\label{prop:nested-cover-criterion}
Retain the preceding hypotheses and suppose that
\begin{equation}\label{eq:nested-cover-condition}
s_\beta\eta(w)\lessdot\eta(w),~\beta\in\Phi^+
\quad\Longrightarrow\quad
M_\beta:=
\langle\mu,\beta\rangle
-\mathbf1_{\{y^{-1}\beta<0\}}=0.
\end{equation}
Then every path $\CP\ne\CP_0$ satisfies $\ell(\eta(w))-\ell(a_{\CP})\geq n_G(\CP)+2.$
Every nonempty stratum $\CS_{\CP} \neq \CS_{\CP_0}$ consequently has dimension
at most $d_w(b)-1$. If $X_{v_{\CP_0}}(b)\neq \emptyset$, then 
\begin{enumerate}[(i)]
    \item $\dim X_w(b)=d_w(b)$;
    \item $\mathscr P(w)_b^{\mathrm{top}}=\{\CP_0\}$, and thus there is a $J_b(F)$-equivariant bijection $\Sigma(X_{v_{\CP_0}}(b))\xrightarrow{\sim}\Sigma^{\text{top}}(X_w(b))$.
\end{enumerate}
\end{proposition}

\begin{proof}
Let $\CP\ne\CP_0$ and put $f=n_G(\CP)\geq 1$.
Recall that $U_m=\eta(w) t^{\gamma_0}$. Suppose the reflection corresponding to the first stationary choice is
$R_i=t^{k\beta^\vee}s_\beta$ for some $1\le k\le M_\beta.$
By Lemma~\ref{lem:stationary-finite}, $s_\beta\eta(w)<\eta(w)$.
Since $M_\beta>0$,
Eq.\ref{eq:nested-cover-condition} implies that this cannot be a covering relation. Since multiplication by a reflection changes
the parity of length, $\ell(\eta(w))-\ell(s_\beta\eta(w))\geq 3$.
Multiplication by each $R_i$ corresponding to subsequent stationary choice decreases the length of the finite Weyl group part
by at least one, while every multiplication choice leaves it
unchanged. Therefore
\[
\ell(\eta(w))-\ell(a_{\CP})
\geq 3+(f-1)=f+2.
\]

For a nonempty stratum, cordiality of the endpoint of its associated path and
\eqref{eq:length-count} give
\[
\begin{aligned}
\dim S_{\CP}(b)
&=\dim X_{v_{\CP}}(b)+\delta_{\CP}\\
&=\frac{
\ell(w)+\ell(a_{\CP})+f
-\langle2\rho,\nu(b)\rangle-\operatorname{def}(b)}2\\
&=d_w(b)
-\frac{\ell(\eta(w))-\ell(a_{\CP})-f}2.
\end{aligned}
\]
Thus every nonempty stratum associated to a path $\CP\neq \CP_0$ has dimension at
most $d_w(b)-1$.

For $\CP_0$, we have $a_{\CP_0}=\eta(w)$ and $n_G(\CP_0)=0$.
If $X_{v_{\CP_0}}(b)\neq\emptyset$, surjectivity of $\pi_{\CP_0}$
makes the ordinary stratum nonempty of dimension $d_w(b)$.
It is therefore the unique top-dimensional stratum.
Since $X_{v_{\CP_0}}(b)$ is equidimensional, Theorem~\ref{thm:main} gives
the asserted equivariant bijection.
\end{proof}
We now use this criterion to give some explicit examples.
\begin{proposition}\label{prop:nested-bijection}
Let $K\subseteq L\subseteq\mathbb S$, and set $y=w_Kw_L$.
Assume $\mu=\rho_{L\setminus K}^\vee\in X_*(T)$.
Then, we have $\mathscr P(w)_b^{\mathrm{top}}=\{\CP_0\}$
whenever $X_{w_Kw_Lt^\gamma}(b)\neq\emptyset$, where $\gamma=\lambda-\mu+w_L\mu$. Thus
\begin{equation*}
 \Sigma\bigl(X_{w_Kw_Lt^\gamma}(b)\bigr)
 \xrightarrow{\sim}
\Sigma^{\text{top}}\bigl(X_{t^\lambda w_Kw_L}(b)\bigr).
\end{equation*}

\end{proposition}
\begin{proof}
It is easy to see that the left descent set of $y=w_Kw_L$ is $L\setminus K$.
Indeed, for $i\in K$, $w_K\alpha_i<0$ and $w_Lw_K\alpha_i>0$;
for $i\in L\setminus K$, $w_K\alpha_i\in\Phi_L^+$ and
$w_Lw_K\alpha_i<0$.  There are no left descents outside $L$.
Moreover $W_K$ fixes $\mu$, so $y^{-1}\mu=w_L\mu$.

Let us verify \eqref{eq:nested-cover-condition}. We note that if $x=1$, the hypothesis simplifies to \[s_\beta y\lessdot y,~\beta\in\Phi^+ \Longrightarrow\langle\mu,\beta\rangle=1.\] Suppose
$s_\beta y\lessdot y$ with $\beta>0$.
The subword property implies $s_\beta y\in W_L$, hence
$\beta\in\Phi_L^+$.  Right multiplication by $w_L$ reverses
Bruhat order in $W_L$, giving $w_K\lessdot s_\beta w_K$.
We then have $w_K=ab, s_\beta w_K=a s_j b$
are reduced expressions, with $a,b\in W_K$ and
$j\in L\setminus K$.  If $j$ were in $K$, the second expression
would give an element of $W_K$ longer than $w_K$.
Thus $s_\beta=a s_j a^{-1}$ and $\beta=a\alpha_j$;
this root is positive because $a\in W_K$ and $j\notin K$.
Since $a^{-1}\mu=\mu$, we obtain $\langle\mu,\beta\rangle
 =\langle\mu,\alpha_j\rangle=1.$
Proposition~\ref{prop:nested-cover-criterion} proves the bijection.
\end{proof}
\begin{corollary}\label{cor:nested-opposite}
Let $K\subseteq L$ and let $K^*$ be the image of $K$ under
the automorphism of $W_L$ defined by
$w_L\alpha_i=-\alpha_{i^*}$.
For $y=w_Lw_K$, the same conclusions as Proposition~\ref{prop:nested-bijection} hold with $J=L\setminus K^*, \mu=\rho_{L\setminus K^*}^\vee,
\gamma=\lambda-\mu+w_L\mu$.

\end{corollary}
\begin{proof}
The identity $w_Lw_K=w_{K^*}w_L$ reduces this to
Proposition~\ref{prop:nested-bijection}.
\end{proof}
For example, let $G=\mathrm{PGL}_4$, $K=\varnothing$,
$L=\mathbb S$, and $\lambda=6\rho^\vee$. Then
$\mu=\rho^\vee$ and $\gamma=\lambda-\mu+w_0\mu=4\rho^\vee$.
The regularity hypothesis holds, and
$X_{w_0t^{4\rho^\vee}}(1)\neq\emptyset$ by the nonemptiness
criterion in Theorem~\ref{prop:xtlambda}. Hence
Proposition~\ref{prop:nested-bijection} gives a $J_1(F)$-equivariant
bijection
\[
\Sigma\bigl(X_{w_0t^{4\rho^\vee}}(1)\bigr)
\xrightarrow{\sim}
\Sigma^{\mathrm{top}}
\bigl(X_{t^{6\rho^\vee}w_0}(1)\bigr).
\]
The other strata need not be empty; they are strictly smaller dimensional.

\printbibliography[heading=bibintoc]

\end{document}